\documentclass[preprint]{elsarticle}
\usepackage{graphicx}
\usepackage{diagbox}
\usepackage{algorithm}
\usepackage{algorithmic}
\usepackage{subcaption}
\usepackage{amsmath}
\usepackage{booktabs}
\usepackage{adjustbox}
\usepackage{fullpage}
\usepackage{amssymb}
\usepackage{hyperref}
\newcounter{bla}

\journal{Journal of Computational and Applied Mathematics}

\begin{document}

\begin{frontmatter}

\title{Spectral-Prior Guided Multistage Physics-Informed Neural Networks for Highly Accurate PDE Solutions}

\author[lab1,lab2]{Yuzhen Li} 
\ead{liyuzhen0201@163.com}
\author[lab1,lab2]{Liang Li\corref{cor1}} 
\ead{plum\_liliang@uestc.edu.cn, plum.liliang@gmail.com}
\author[lab3]{St\'ephane Lanteri} 
\ead{stephane.lanteri@inria.fr}
\author[lab2]{Bin Li}
\ead{libin@uestc.edu.cn}

\address[lab1]{School of Mathematical Sciences, University of
               Electronic Science and Technology of China, 611731,
               Chengdu, P.R. China}
               
\address[lab2]{Key Laboratory of Large-Scale Electromagnetic Industrial Software, Ministry of Education, , 611731,
               Chengdu, P.R. China}

\address[lab3]{Universit\'e C\^ote d'Azur, Inria, CNRS, LJAD, Sophia Antipolis, France}

\cortext[cor1]{Corresponding author}

\begin{abstract}
Physics-Informed Neural Networks (PINNs)  are becoming a popular method for solving PDEs, due to their mesh-free nature and their ability to handle high-dimensional problems where traditional numerical solvers often struggle. Despite their promise, the practical application of PINNs is still constrained by several factors, a primary one being their often-limited accuracy.  This paper is dedicated to enhancing the accuracy of PINNs by introducing spectral-prior guided multistage strategy. We propose two methods: Spectrum-Informed Multistage Physics-Informed Neural Networks (SI-MSPINNs) and Multistage Physics-Informed Neural Networks with Spectrum Weighted Random Fourier Features (RFF-MSPINNs). The SI-MSPINNs integrate the core mechanism of Spectrum-Informed Multistage Neural Network (SI-MSNNs) and PINNs, in which we extract the Dominant Spectral Pattern (DSP) of residuals by the discrete Fourier transform. This DSP guides the network initialization to alleviate spectral bias, and gradually optimizes the resolution accuracy using a multistage strategy. The RFF-MSPINNs combines random Fourier features with spectral weighting methods, dynamically adjusting the frequency sampling distribution based on the residual power spectral density, allowing the network to prioritize learning high-energy physical modes. Through experimental verification of the Burgers equation and the Helmholtz equation, we show that both models significantly improve the accuracy of the original PINNs.
   \end{abstract}
   
   
\end{frontmatter}

\section{Introduction}
\label{Introduction}
In recent years, the rapid development of deep learning has brought new opportunities to scientific computing. In 2019, Raissi et al. proposed Physics-Informed Neural Networks (PINNs) \cite{Raissi}, which embed physical laws into the neural network training process, improving the interpretability of models and solving PDEs and inverse problems without relying on large amounts of data. This mesh-free approach allows them to handle complex geometries and high-dimensional spaces with greater flexibility than traditional methods like the finite element method (FEM) \cite{Raissi,Gomez-Serrano}.  By combining data-driven methods with physical constraints \cite{chen2021physics,FU2023115771}, PINNs provide powerful modeling capabilities, particularly suitable for solving high-dimensional nonlinear systems. Chen et al. \cite{chen2021physics} proposed a physics-informed machine learning framework for the reduced-order modeling of parametric PDEs, demonstrating its effectiveness through steady and unstable problems. 

However, the widespread adoption of PINNs is still hindered by significant limitations, with sub-optimal accuracy being a primary concern \cite{Brockherde,Westermayr,Chattopadhyay}. This issue arises mainly from the complex and often problematic optimization dynamics of the training process \cite{Karniadakis}. The composite loss function, which combines terms for PDE residuals, boundary, and initial conditions, frequently suffers from imbalanced gradients. Furthermore, the inherent spectral bias of neural networks makes them slow to learn high-frequency or multi-scale features of a solution, a critical weakness for many complex physical systems \cite{Rahaman,Xu,Jacot}. These optimization challenges make PINNs difficult to have high accuracy, and thus less robust and reliable than mature numerical methods like FEM or FVM for many applications \cite{Mojgani,Hassanzadeh,Perdikaris,Udrescu}.

In order to improve the accuracy of the Neural Network (NN) based methods, Ng et al. \cite{MSNN} proposed the Multistage Neural Networks (MSNNs) method, which divides the neural network training process into multiple stages. Each stage learns the residual of the previous stage (i.e. the difference between the objective function and the output of the trained network) through a new neural network. And by introducing scale factors to adjust weights and other methods, the network settings are optimized for the high-frequency characteristics of residuals to alleviate the spectral bias \cite{ZhangJ} of the neural network and achieve a high-precision approximation of the objective function. However, the convergence efficiency decreases significantly in high-dimensional problems. The author applied it to regression problems and ordinary/partial differential equations to verify its effectiveness. Ng et al. \cite{SI-MSNN} proposed an improved Spectrum-Informed Multistage Neural Network (SI-MSNN) method based on MSNN. This method further alleviates the spectral bias of the neural network through a spectral information guided initialization strategy, enabling the network to quickly converge to the spectral model of the objective function. By combining the strategy of learning the residuals of the previous stage in multi-stage training, a high-precision approximation of the two-dimensional regression problem was ultimately achieved, but it was not applied to the PINNs framework to solve the PDEs problem.

Based on the above considerations, this paper integrates the SI-MSNN method into the PINNs framework and establishes the Spectrum-Informed Multistage Physics-Informed Neural Networks (SI-MSPINNs). SI-MSPINNs utilizes a spectrum guided initialization strategy to enable the network to more accurately capture high-frequency features when learning physics equation residuals, thereby improving the accuracy and efficiency of solving complex physics problems. In addition, a spectral weighting method incorporating random Fourier features (RFF) was introduced to initialize MSNN, and an Multistage Physics-Informed Neural Networks with Spectrum Weighted Random Fourier Features (RFF-MSPINNs) model was established. By assigning spectral weights to random Fourier features, the accuracy of MSNN in handling 2D problems in PINNs has been improved, and the spectral bias problem in neural networks has also been solved.

The remainder of this paper is organized as follows. Section \ref{sec:RW} briefly introduces the background and theoretical knowledge of related works such as PINNs and MSNNs. Section \ref{sec:SI-MSPINNs} introduces the framework and algorithm of SI-MSPINNs. Section \ref{sec:RFF-MSPINNs} introduces the framework and algorithm of RFF-MSPINNs. Section \ref{sec:Experiments} applies the two newly proposed models to the Burgers equation and the Helmholtz equation, respectively. Section \ref{sec:Conclusion} provides a summary of the entire text and offers prospects for future work.

\section{Related Work}\label{sec:RW}
\subsection{PINNs: Physics-Informed Neural Networks}
PINNs is a PDEs solving paradigm that deeply integrates physical laws into deep learning frameworks. Its core idea is to achieve unified modeling of physical laws and data constraints through a dual loss collaborative mechanism. This method constructs an approximate solution function for neural networks and synchronously embeds control equation residuals and boundary/initial condition constraints in the loss function. And by using automatic differentiation techniques to calculate higher-order derivative terms, without relying on grid discretization in traditional numerical methods, the overall loss is minimized by optimizing network parameters, ultimately obtaining a solution that conforms to both data and physical laws.

Consider the following partial differential equation defined on a spatial domain $\Omega \subset \mathbb{R}^d$ with boundary $\partial\Omega$:
\begin{equation}  
\begin{aligned}  
&\mathcal{L}\bigl(u(\mathbf{x})\bigr) = f(\mathbf{x}) & \quad & \textbf{x} \in \Omega, \quad t \in [0,T], \\
&\mathcal{B}(u(\textbf{x},t);\mu )= 0 & \quad & \textbf{x} \in \partial \Omega, \quad t \in [0,T], \\
&\mathcal{I}(u(\textbf{x},0)) = u_0(\textbf{x};\mu) & \quad & \textbf{x} \in \Omega
\end{aligned}  
\end{equation} 
where $\mathcal{L}$ is a (nonlinear) differential operator, $\mathcal{B}$ and $\mathcal{I}$ are operators for boundary condition and initial condition, respectively. The goal of PINNs is to approximate the solution of the physical equation using a neural network $P^{\theta}(\mathbf{x},t)$ where $\theta$ represents the parameters of the neural network, including weights and bias vectors. A loss function is defined to learn the parameters $\theta$ by minimizing the value of the loss function. The loss function in PINNs is defined as a weighted sum of the residuals of the physical equation, the boundary conditions and the initial conditions \cite{Raissi}:
\begin{equation}  
\begin{aligned}
    \mathcal{L}_{\text{PINN}}(P^{\theta}) &= \frac{\lambda_f}{|C_f|} \sum_{(\textbf{x},t) \in C_f} \left\| \mathcal{L}\bigl(u(\mathbf{x})\bigr) + f(\mathbf{x}) \right\|^2_2 +\\ &\quad\frac{\lambda_b}{|C_\partial|} \sum_{(\textbf{x},t) \in C_\partial} \left\| \mathcal{B}(P^{\theta})(\textbf{x},t) \right\|^2_2 + \frac{\lambda_i}{|C_i|} \sum_{\textbf{x} \in C_i} \left\| P^{\theta}(\textbf{x}, 0) - u_0(\textbf{x}) \right\|^2_2,
\end{aligned}
\end{equation} 
where $C_f \subset \Omega \times [0,T]$, $C_\partial \subset \partial \Omega \times [0,T]$ and $C_i \subset \Omega$ represent the sampling spaces, $|C_f|$, $|C_\partial|$ and $|C_i|$ denote the number of points of insertion for the equation, boundary, and initial conditions, respectively, and $\lambda_f$, $\lambda_b$ and $\lambda_i$ represent the weights of the three loss terms. In this way, PINNs combine the data-driven approach with the physics-based approach to improve the accuracy and reliability of the neural network model. 

\subsection{MSNNs: Multistage Neural Networks}
MSNNs effectively alleviate the spectral bias of traditional neural networks and significantly improve the accuracy of the approximation of functions by dividing the training process into multiple stages \cite{Wang,LiuX}, each stage using a new network to fit the residual of the previous stage. The key idea is to normalize the residual magnitude and adjust the scale factor based on the residual characteristics of different stages, to better initialize the network. Below, the processes of these two core mechanisms are briefly introduced, respectively.

\subsubsection{Magnitude Normalization}
When traditional neural networks fit residuals, if the residual magnitude is much smaller than the original data, conventional weight initialization methods such as Xavier are difficult to effectively capture small signals, resulting in limited accuracy in subsequent network fitting \cite{MSNN}. To address this issue, MSNN adopts the following strategy.
\begin{itemize}
  \item Residual normalization: Calculate the root mean square value (RMS) of the residual $e_1(x)$ from the previous stage, defined as
\begin{equation}
\epsilon_1 = \text{RMS}(e_1(x)) = \sqrt{ \frac{1}{N_d} \sum_{i=1}^{N_d} \left[ e_1(x^{(i)}) \right]^2 } = \sqrt{ \frac{1}{N_d} \sum_{i=1}^{N_d} \left[ u_g^{(i)} - u_0(x^{(i)}) \right]^2 },
\end{equation}
where $u_0$ represents the neural network trained in the first stage on data points $N_d$ and $u_g$ represents the objective function. Then normalize the residual to $e_1(x)/\epsilon_1$ as the training target for the second stage network.
  
  \item Scale integration: The output of the second stage network $u_1(x)$ must be multiplied by a normalization coefficient of $\epsilon_1$, and then combined with the previous stage network $u_0(x)$, which is $u_c^{(1)}(x) = u_0(x) + \epsilon_1 u_1(x)$. This process can be recursively extended to more stages, and the final model is $u_c^n(x) = \sum_{j=0}^n \epsilon_j u_j(x)$, where $\epsilon_j$ is the magnitude coefficient of the residuals in each stage.
\end{itemize}

\subsubsection{Frequency Adaptation}
\begin{figure}[htbp]
\centering
\includegraphics[scale=0.8]{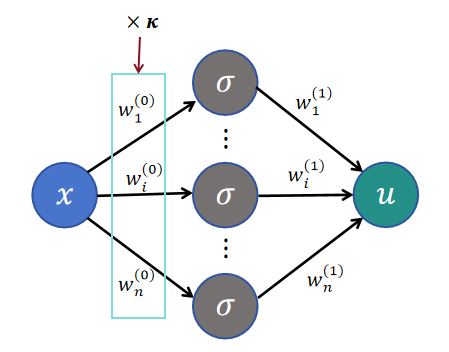}
\caption{Schematic diagram of scale factor $\kappa$ applied to a single hidden layer neural network.}
\label{MSNN}
\end{figure}

Traditional neural networks suffer from spectral bias, making it difficult to fit high-frequency features. In multistage training, the high-frequency components of residuals are often more significant. Therefore, for high-frequency residuals, it is necessary to increase the weight scale from the input layer to the first hidden layer. MSNN accelerates the convergence of weights to the optimal value required for high frequencies (approximately $O(2\pi f_d)$, where $f_d$ is the residual dominant frequency) by introducing scale factor $\kappa$ and multiplying the weights of the layer by $\kappa$. The schematic diagram is shown in Figure \ref{MSNN}, considering a shallow neural network with a single input, a single output, and a hidden layer, where $w_i^{(0)}$ represents the weight between the input layer and the hidden layer, and $w_i^{(1)}$ corresponds to the weight between the hidden layer and the output layer, and $\sigma$ represents the activation function. Multiply the weights $w_i^{(0)}$ between the input layer and the first hidden layer by the scale factor $\kappa$.

\subsection{SI-MSNNs: Spectrum-Informed Multistage Neural Networks \cite{SI-MSNN}}
SI-MSNNs is an improved model for high-precision function approximation problems. This method aims to solve the problem of difficulty in capturing high-frequency features and limited accuracy of MSNNs in high-dimensional scenes. By integrating spectral information of the objective function to optimize network initialization, the accuracy limit of high-dimensional problems is significantly improved. The core innovation of SI-MSNNs lies in the spectrum-informed initialization strategy, which extracts the dominant frequency features of the objective function based on discrete Fourier transform (DFT) and embeds them into the first layer parameters of the network. The specific steps are as follows:
\begin{enumerate}
  \item Calculate the discrete Fourier transform $\tilde{u}(k)$ of the objective function $u$ or residual, and obtain its polar coordinate form $\tilde{u}(k) = a(k) e^{i\theta(k)}$, where $a(k)$ is the amplitude and $\theta(k)$ is the phase.
  \item Select the top $n_f$ Fourier modes with the largest amplitude ($k^{(j)}$), and calculate the normalized amplitude $a^{(j)}$ and the corresponding phase $\theta^{(j)}$. For two-dimensional problems, use the conjugate symmetry of Fourier transform to screen patterns with $k_y>0$ to avoid redundancy.
  \item Map the selected patterns to the first layer of the neural network, the Fourier feature embedding layer, which has an activation function of $\gamma(x) = A \cos(Bx + b)$ and directly encodes the main spectral information of the objective function, where the weight matrix $B = [k^{(1)}, \cdots, k^{(n_f)}]^{\mathrm{T}}$ corresponds to the dominant frequency, the bias vector $b = [\theta^{(1)}, \cdots, \theta^{(n_f)}]^{\mathrm{T}}$ corresponds to the phase, and the amplitude parameter $A = [\alpha^{(1)}, \cdots, \alpha^{(n_f)}]^{\mathrm{T}}$ corresponds to the feature strength.
\end{enumerate}

In summary, SI-MSNNs provide an efficient framework for high-precision function approximation and validate their effectiveness in multiscale scientific problems.

\section{Spectrum-Informed Multistage Physics-Informed Neural Network}\label{sec:SI-MSPINNs}
This section integrates the core ideas of SI-MSNNs with the PINNs framework to propose the SI-MSPINNs model, which solves the accuracy bottleneck of traditional PINNs in solving high-frequency PDEs. This model inherits the residual learning mechanism of MSNN and uses the SI-MSNNs method to extract the spectral features of residuals through discrete Fourier transform, guiding the initialization of network weights and effectively alleviating spectral bias. At the same time, the physical constraint core of PINNs is retained, and the residual of partial differential equations is incorporated into the loss function to achieve collaborative optimization of data and physical laws.

Traditional PINN simultaneously fits data and physical equation residuals through a single-stage network, which is susceptible to spectral bias limitations and difficult to capture high-frequency physical features. SI-MSPINNs borrows the multistage training logic of SI-MSNN and decomposes PDE solving into a progressive residual learning process.
Given PDE:
\begin{equation}
\mathcal{L}\bigl(u(\mathbf{x})\bigr) = f(\mathbf{x}), \quad \mathbf{x} \in \Omega,
\label{eq:pde}
\end{equation}
and boundary condition:
\begin{equation}
\mathcal{B}\bigl(u(\mathbf{x})\bigr) = g(\mathbf{x}), \quad \mathbf{x} \in \partial\Omega.
\label{eq:bc}
\end{equation}
First, train the basic PINN network $u_0$, where the loss function includes data fitting loss (initial/boundary conditions) and PDE residual loss, consistent with traditional PINN,
\begin{equation}
\mathcal{L}_\text{PINN} = \|\mathcal{L}(u_0) - f\|^2 + \|\mathcal{B}(u_0) - g\|^2.
\end{equation}
Then, calculate initial residual  \( r_0 = \mathcal{L}(u_0) - f \) Stage $n \geq 1$. Specifically, perform DFT on residual $r_{n-1}$ to extract the first $n_f$ main modes (amplitude $a_j$ and corresponding phase $\theta_j$). Then build a spectral embedding layer:
\begin{equation}
\gamma(\mathbf{x}) = [\alpha_j \cos(\mathbf{k}_j \cdot \mathbf{x} + \theta_j)]^{\mathrm{T}}, \quad j = 1, \ldots, n_f.
\end{equation}
Finally, calculate the RMS value $\epsilon_{n-1}$ of the residual and train $u_n$ to fit the normalized residual $r_{n-1}/\epsilon_{n-1}$ to $\min \left\| \mathcal{L}(u_n) - \frac{r_{n-1}}{\|\epsilon_{n-1}\|} \right\|^2$. Finally update total solution $u_{\text{total}} = u_{\text{total}} + \epsilon_{n-1} u_{\theta_n}$ and residual $\epsilon_k = \mathcal{L}(u_{\text{total}}) - f$.
The algorithm and flowchart of SI-MSPINNs are given by Algorithm \ref{alg:SI-MSPINNs} and Figure \ref{fig:SIflowchart} respectively.

\begin{algorithm}[H]
\caption{Training Scheme for SI-MSPINNs}
\label{alg:SI-MSPINNs}
\begin{algorithmic}[1]
\REQUIRE PDE operator $\mathcal{L}[\cdot]$, training data $\{(x_i, f(x_i))\}_{i=1}^{N_d}$
\STATE Train the 0-th stage PINN network $u_0$ with loss $\mathcal{L}_{\text{PINN}}(u_0)$, yielding the initial solution $u_0$.
\STATE Compute the initial residual $r_0 = u - u_0$ and the normalized coefficient $\epsilon_0 = \sqrt{\dfrac{1}{N_d} \sum_{i=1}^{N_d} r_{0,i}^2}$.

\WHILE{$n \leq s$}
    \STATE \textbf{Spectrum Extraction:} Perform DFT on ${r}_{n-1}$ to obtain the spectrum $\tilde{r}_{n-1}(k) = a(k) e^{i\theta(k)}$, and select the top $n_f$ dominant spectral modes $\{k^{(j)}, \alpha^{(j)}, \theta^{(j)}\}_{j=1}^{n_f}$.
    \STATE \textbf{Network Initialization:} Initialize the first-layer weights $B$, biases $b$, and amplitudes $A$ of the $n$-th stage network $u_n$ using $\{k^{(j)}, \alpha^{(j)}, \theta^{(j)}\}$; other layers are initialized via Xavier initialization.
    \STATE \textbf{Train $u_n$:} Minimize the loss $\mathcal{L}_{\text{PINN}}(u_0)$ to obtain $u_n$.
    \STATE \textbf{Update Residual:}  $r_n(x) = \mathcal{L}[u_n](x) - f(x)$, and compute $\epsilon_n = \sqrt{\dfrac{1}{N_d} \sum_{i=1}^{N_d} r_{n,i}^2}$.
    \STATE Set $n = n + 1$.
\ENDWHILE

\ENSURE The final solution: $u_s = \sum_{i=0}^s \epsilon_i u_i$.
\end{algorithmic}
\end{algorithm}

\begin{figure}[H]
\centering
\includegraphics[scale=0.3]{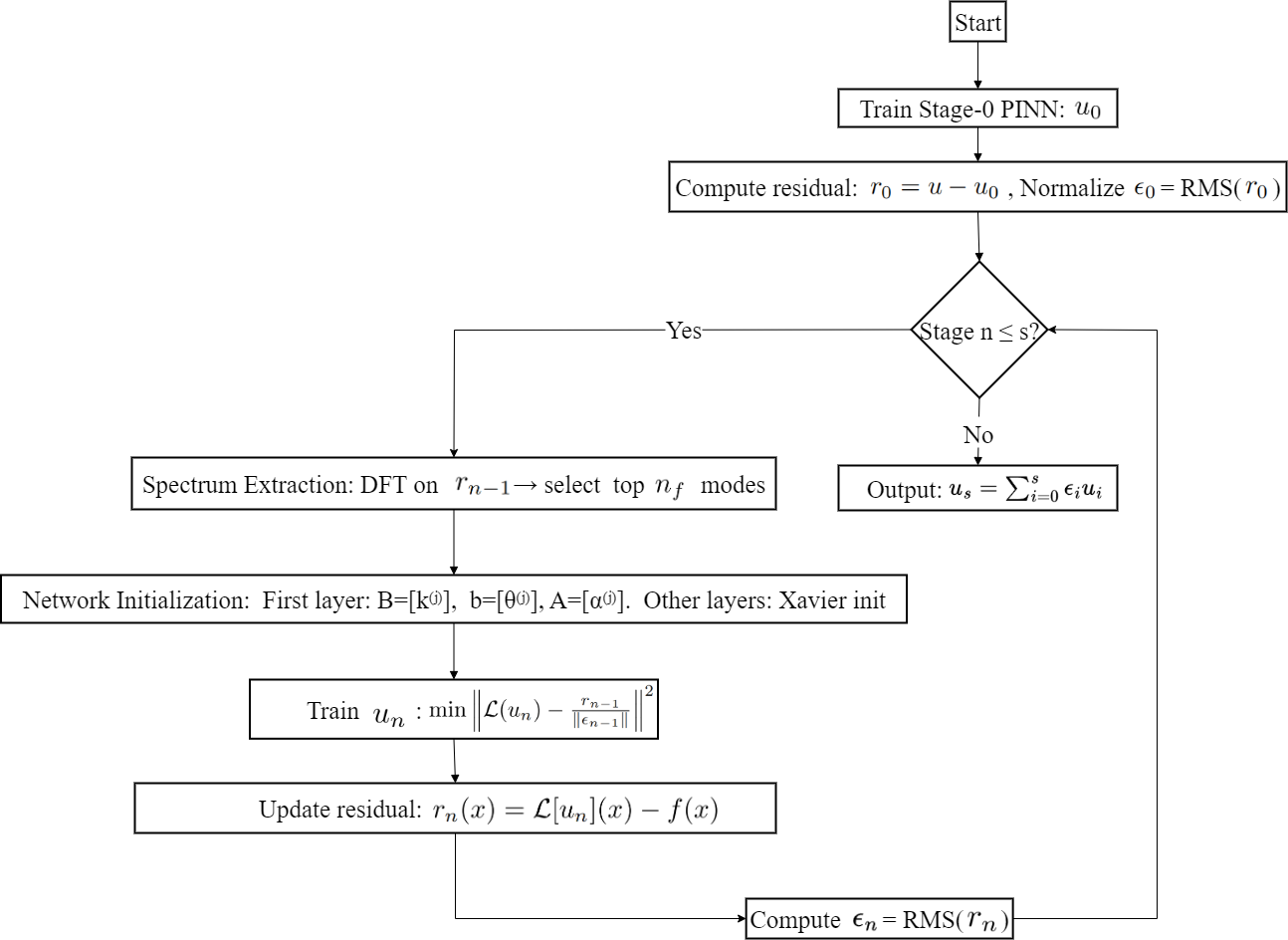}
\caption{The flowchart of the SI-MSPINNs.}
\label{ig:SIflowchart}
\end{figure}

The key idea of the proposed SI-MSPINNs method in this section is to integrate spectrum information initialization and multistage physical constraints. 

\section{Multistage Physics-Informed Neural Network with Spectrum Weighted Random Fourier Features}\label{sec:RFF-MSPINNs}
RFF-MSPINNs integrates the spectral weighting method of random Fourier features (RFF) with MSNNs into the PINNs framework, significantly improving the accuracy and efficiency of PINNs. RFF-MSPINNs utilizes the power spectral density (PSD) of objective function residuals to guide RFF frequency sampling, enabling the neural network to prioritize learning the energy dominated physical modes and break through the bottleneck of high-frequency learning.

RFF is a technique that constructs high-dimensional feature maps using sine and cosine functions to approximate the kernel functions by randomly sampling frequency vectors and phases \cite{Rahimi,ChenW}. Its core idea is to map the low dimensional inputs to high-dimensional space, enhance the model's ability to learn high-frequency features, and alleviate the spectral bias of neural networks. RFF uses the Equation \ref{rff} for random frequency sampling to approximate shift invariant kernels,
\begin{equation}
\phi(\mathbf{x}) = \sqrt{\frac{2}{m} \left[ \cos\left( \boldsymbol{\omega}_1^T \mathbf{x} + b_1 \right), \ldots, \cos\left( \boldsymbol{\omega}_m^T \mathbf{x} + b_m \right) \right]^T}
\label{rff}
\end{equation}
where  \( \boldsymbol{\omega}_j \sim \mathcal{N}(0, \sigma^2) \), \( b_j \sim \mathcal{U}[0, 2\pi] \). Insufficient high-frequency sampling in traditional RFF may result in spectral bias.  

In order to address this, PSD \cite{Li,WangY,Chew} is introduced, which is a physical quantity that describes the energy distribution of a signal or function in the frequency domain. It quantitatively characterizes the contribution of different frequency components to the total energy of the signal or function. The larger the PSD value at a certain frequency, the higher the proportion of that frequency component in the signal or function \cite{Shanker}. It is a key indicator for analyzing the multiscale and spectral characteristics of signals.

In the proposed RFF-MSPINNs, we adopt a spectrum weighted RFF mechanism, whose core idea is to adjust the sampling distribution based on the residual PSD and prioritize the energy dominant frequency. The key process consists of the following four steps:
\begin{enumerate}
  \item Firstly, perform DFT on the residual function $r_n(x)$, that is, $\tilde{r}_n(k) = \mathcal{F}\{r_n(x)\}$. And calculate the residual PSD $P(k) = \left| \tilde{r}_n(k) \right|^2$.
  \item Then use Equation \ref{psd} to normalize PSD to a probability distribution, in order to construct a frequency sampling distribution
  \begin{equation}
p(k) = \frac{P(k)}{\sum_k P(k)}.
\label{psd}
\end{equation}
  \item Next, in order to utilize the spectral information of the residuals, we change the frequency sampling distribution to be proportional to the PSD of the residuals. As shown in Equation \ref{psd1}, extract $m$ frequencies from $p(k)$.
  \begin{equation}
k_j \sim p(k),\ j = 1, 2, \ldots, m
\label{psd1}
\end{equation}
\item Finally, initialize the RFF layer by constructing frequency matrix $\mathbf{B} \in \mathbb{R}^{m \times d}$:
\begin{equation}
\mathbf{B} = \begin{bmatrix} k_1 \\ k_2 \\ \vdots \\ k_m \end{bmatrix}.
\end{equation}
Randomly generate phase shifts $b_j \sim \mathcal{U}[0, 2\pi]$, and finally defining the feature mapping:
\begin{equation}
\gamma(x) = \cos(2\pi \mathbf{B} x + \mathbf{b}).
\end{equation}
\end{enumerate}

\begin{figure}[H]
\centering
\includegraphics[scale=0.3]{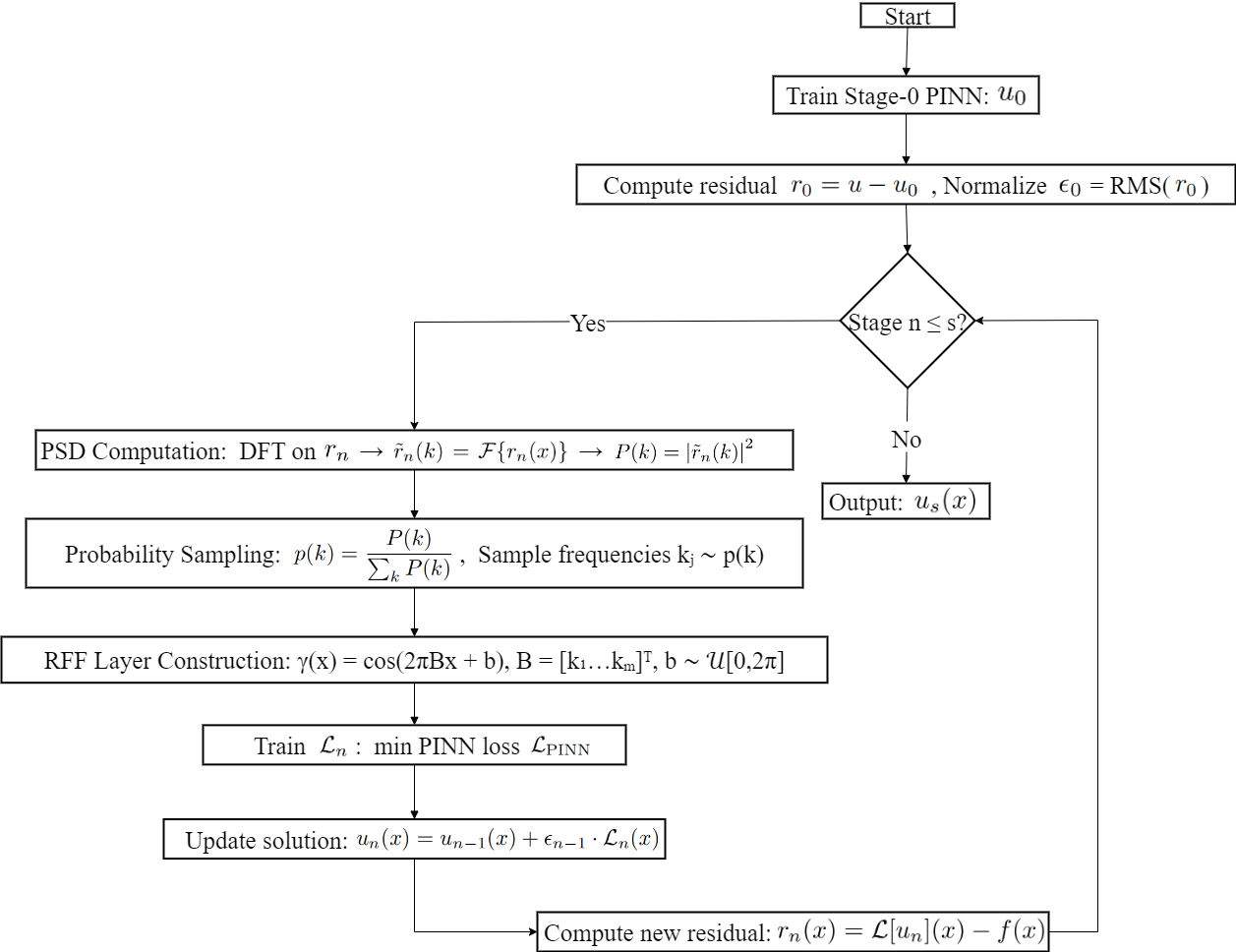}
\caption{The flowchart of RFF-MSPINNs.}
\label{RFFflowchart}
\end{figure}

This RFF module enhances the network's ability to learn high-frequency components by using a frequency sampling distribution proportional to the residual PSD. By sampling frequency based on power spectral density, we ensure that more features are allocated on important frequencies (high-energy frequencies) and fewer features are allocated on low-energy frequencies. In this way, the neural network has features that match the frequency spectrum of the objective function during initialization, allowing it to learn the main frequency components of the objective function faster, especially those with higher energy.
The algorithm and flowchart of SI-MSPINNs are given by Algorithm \ref{RFF-MSPINNs} and Figure \ref{RFFflowchart} respectively.

\begin{algorithm}[H]  
\caption{Training Scheme for RFF-MSPINNs}  
\label{RFF-MSPINNs}  
\begin{algorithmic}[1]  

    \REQUIRE PDE operator $\mathcal{L}[\cdot]$, training data $\{(x_i, f(x_i))\}_{i=1}^{N_d}$  
    \STATE \textbf{Initialization}  
    \STATE \hspace{\algorithmicindent} Train the 0-th stage PINN network $u_0$ with loss $\mathcal{L}_{\text{PINN}}(u_0)$, yielding the initial solution $u_0$.  
    \STATE \hspace{\algorithmicindent} Compute the initial residual $r_0 = u - u_0$ and the normalized coefficient $\epsilon_0 = \sqrt{\dfrac{1}{N_d} \sum_{i=1}^{N_d} r_{0,i}^2}$.

    \WHILE{$n \leq s$}  
        \STATE Compute the PSD $P_{n-1}(k)$ of $r_{n-1}(x)$ via DFT.  
        \STATE Sample RFF frequencies $\{\mathbf{b}_i\}_{i=1}^M$ according to $p(k) \propto P_{n-1}(k)$, and phases $b_i \sim \mathcal{U}(0, 2\pi)$.  
        \STATE Construct network $\mathcal{N}_n$ containing the RFF layer $\gamma(x)$ and $L$ hidden layers.  
        \STATE Normalize the residual: $\hat{r}_{n-1}(x) = \dfrac{r_{n-1}(x)}{\epsilon_{n-1}}$, where $\epsilon_{n-1} = \text{RMS}(r_{n-1})$.  
        \STATE Train $\mathcal{L}_n$ by minimizing the PINN loss $\mathcal{L}_{\text{PINN}}(u_n)$.  
        \STATE Update the solution: $u_n(x) = u_{n-1}(x) + \epsilon_{n-1} \cdot \mathcal{L}_n(x)$.  
        \STATE Compute the new residual: $r_n(x) = \mathcal{L}[u_n](x) - f(x)$.  
        \STATE Set $n = n + 1$.
\ENDWHILE  

    \ENSURE Final solution $u_s(x)$  

\end{algorithmic}  
\end{algorithm}  

The SI-MSNN directly select the first $n_f$ maximum amplitude frequencies and explicitly sets their amplitude and phase, which is equivalent to constructing a deterministic initialization using these frequency components. While, the RFF-MSPINNs select frequencies through probability sampling, and the probability of each frequency being selected is proportional to its energy. The initialization constructed in this way contains randomness, but statistically tends to favor important frequencies. It should be noted that RFF-MSPINNs guides the network to focus on important frequencies while maintaining randomness, which is a method that falls between completely deterministic (SI-MSNN) and completely random (traditional RFF). It may be more scalable in high-dimensional problems because it does not require sorting of all frequencies (SI-MSNN requires sorting) but only sampling based on probability distribution.

\section{Experiments}\label{sec:Experiments}
This section applies the SI-MSPINNs and RFF-MSPINNs proposed above to the Burgers equation and the Helmholtz equation, respectively, to verify the performance of SI-MSPINNs and RFF-MSPINNs. All experiments are implemented using Torch 2.5.1 on a workstation with an Intel i9-10900 CPU and an NVIDIA Quadro RTX 5000 GPU. Codes and the corresponding data are available on GitHub \footnote{available at https://github.com/liyuzhen0201/MPINN}.

\subsection{Burgers Equation}
Consider the Burgers equation:
\begin{equation}
\frac{\partial u}{\partial t} + u \frac{\partial u}{\partial x} = \frac{\partial^2 u}{\partial x^2}, \quad (x,t) \in [-1,1] \times [0,1],
\label{eq:pde_example}
\end{equation}
with Dirichlet boundary conditions and initial conditions $u(-1,t) = u(1,t) = 0$ and $u(x,0) = -\sin(\pi x)$.

Set the number of training residual point samples within the domain to 2540, the number of boundary sampling points to 80, and the number of initial condition residual points. We use a three-stage architecture. The three stages of this case all use fully connected neural networks with a depth of 4 (equivalent to 3 hidden layers) and a width of 20. Adam+L-BFGS is used as the optimizer for each of the three stages. 

\subsubsection{Results of the SI-MSPINNs}
\begin{figure}[ht]  
\centering
\begin{minipage}{\textwidth}
\centering
\footnotesize (a) First-stage

\begin{subfigure}{0.32\textwidth}
    \centering
    \includegraphics[width=\textwidth]{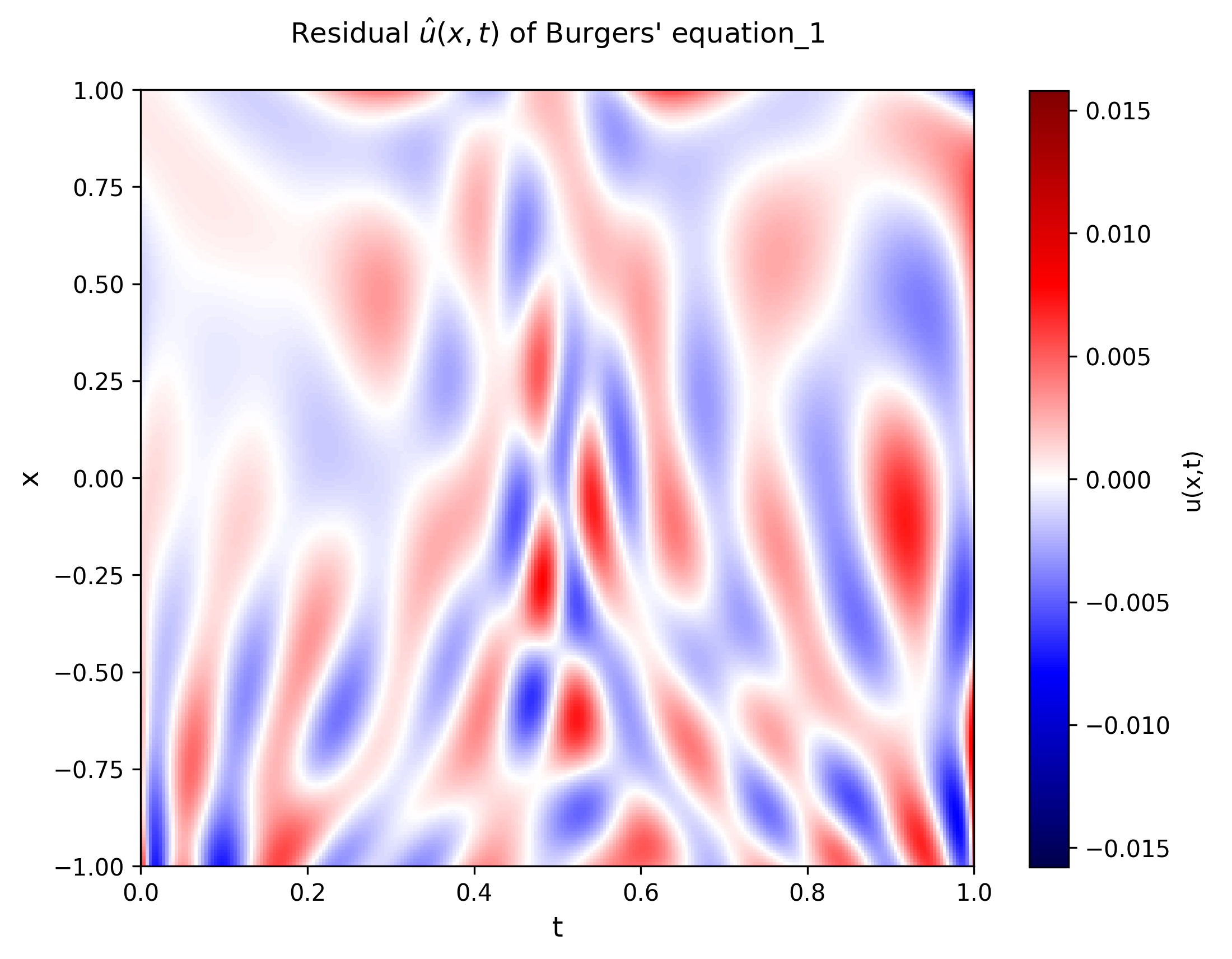}
\end{subfigure}
\hfill
\begin{subfigure}{0.32\textwidth}
    \centering
    \includegraphics[width=\textwidth]{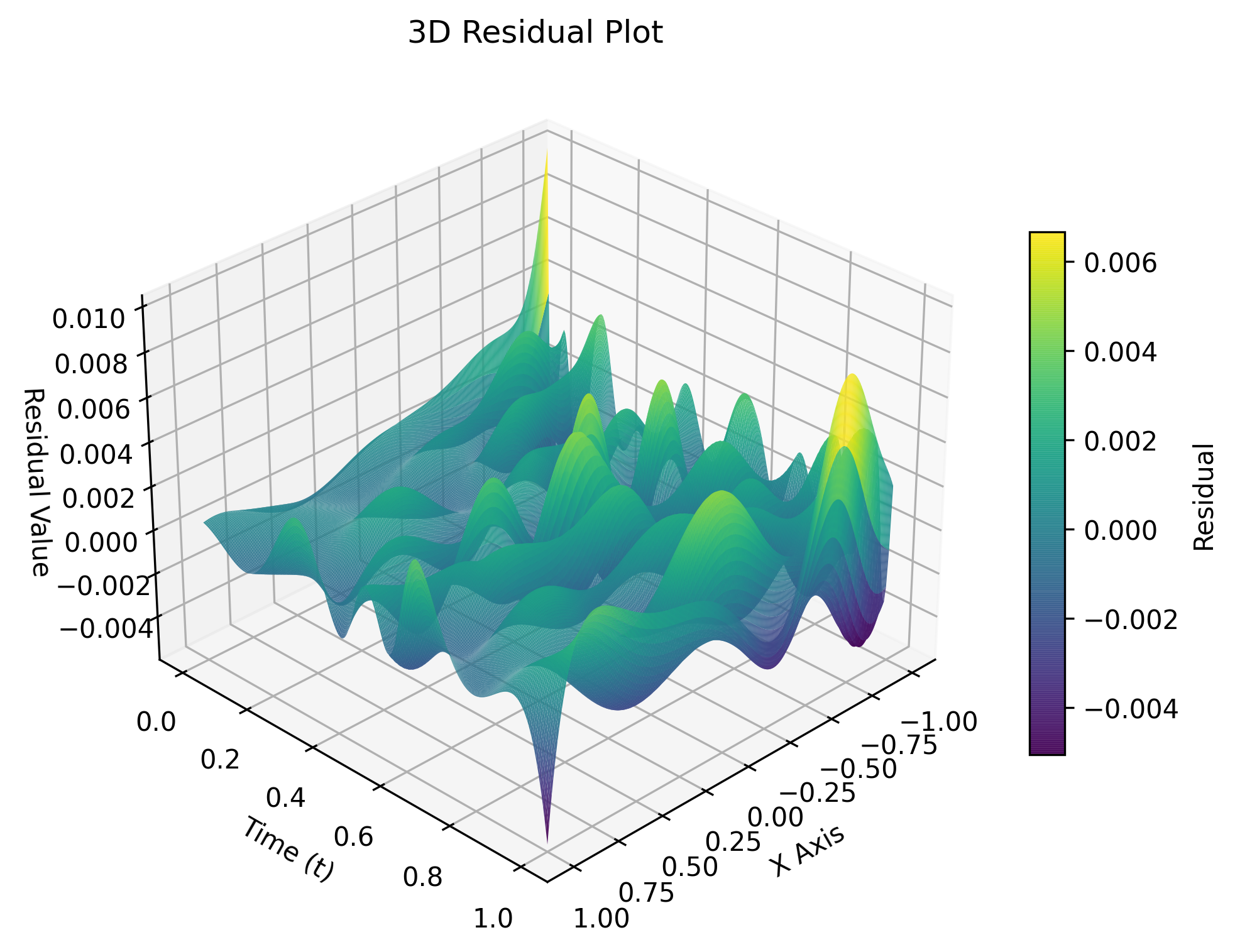}
\end{subfigure}
\hfill
\begin{subfigure}{0.32\textwidth}
    \centering
    \includegraphics[width=\textwidth]{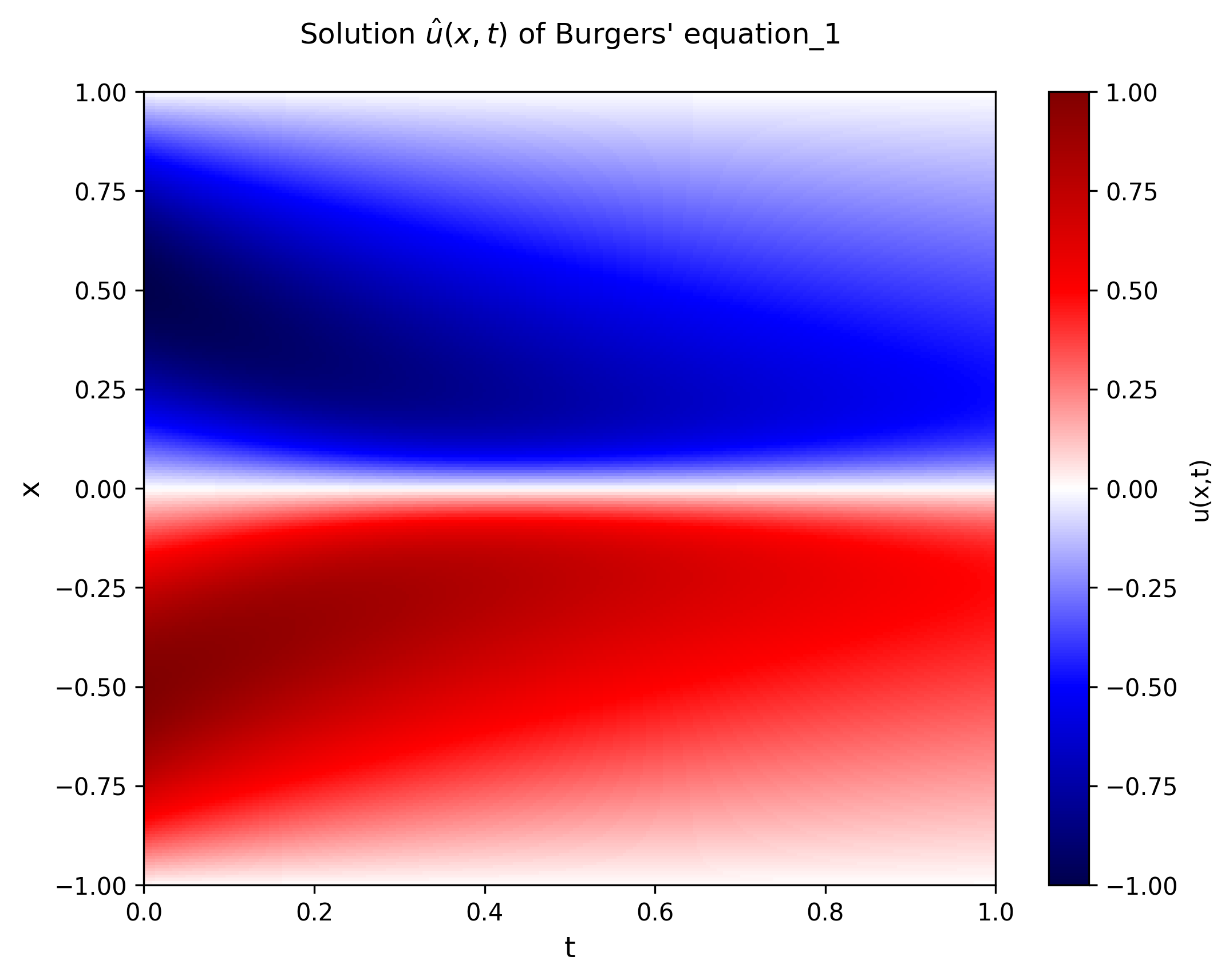}
\end{subfigure}
\end{minipage}

\vspace{1em} 

\begin{minipage}{\textwidth}
\centering
\footnotesize (b) Second-stage

\begin{subfigure}{0.32\textwidth}
    \centering
    \includegraphics[width=\textwidth]{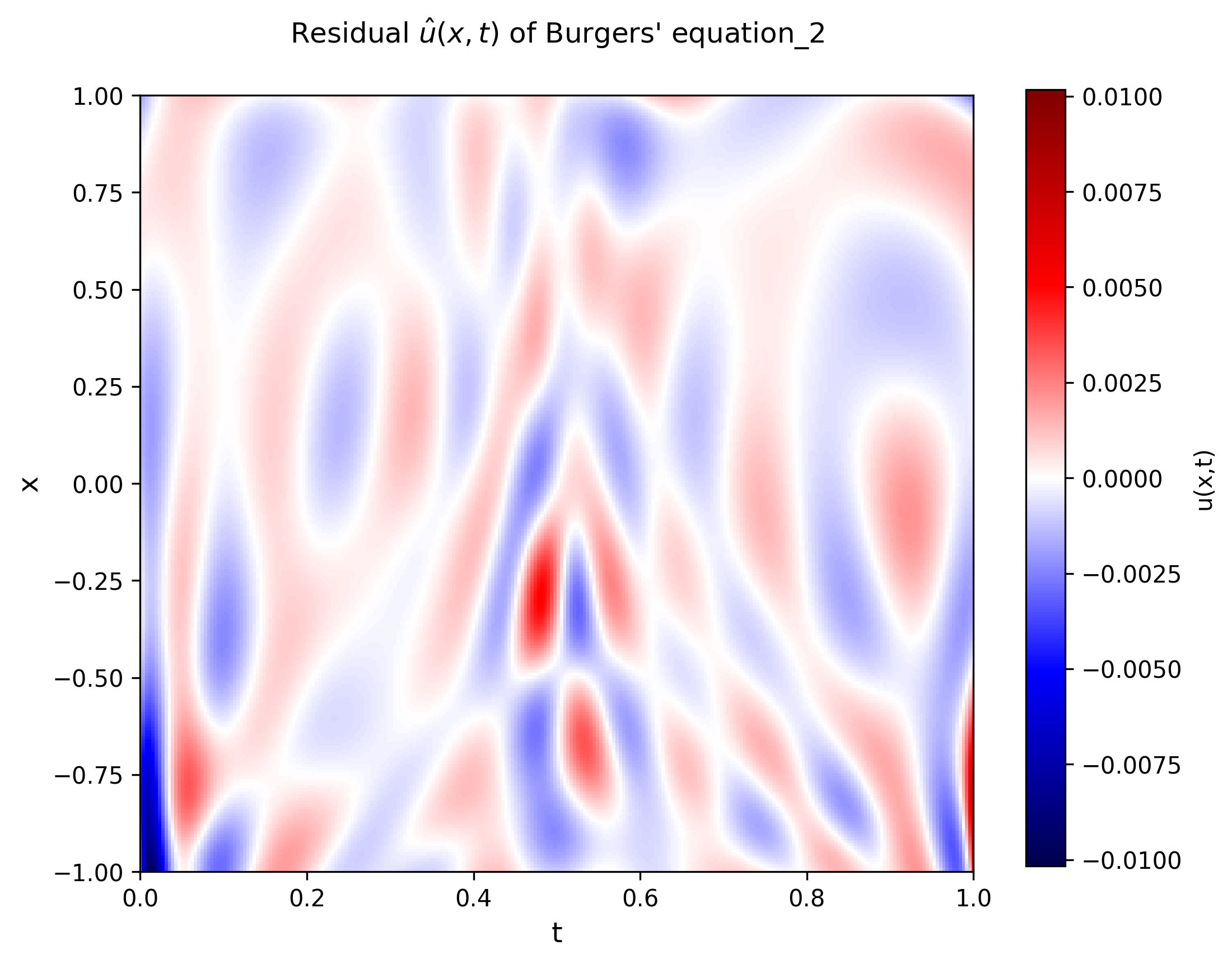}
\end{subfigure}
\hfill
\begin{subfigure}{0.32\textwidth}
    \centering
    \includegraphics[width=\textwidth]{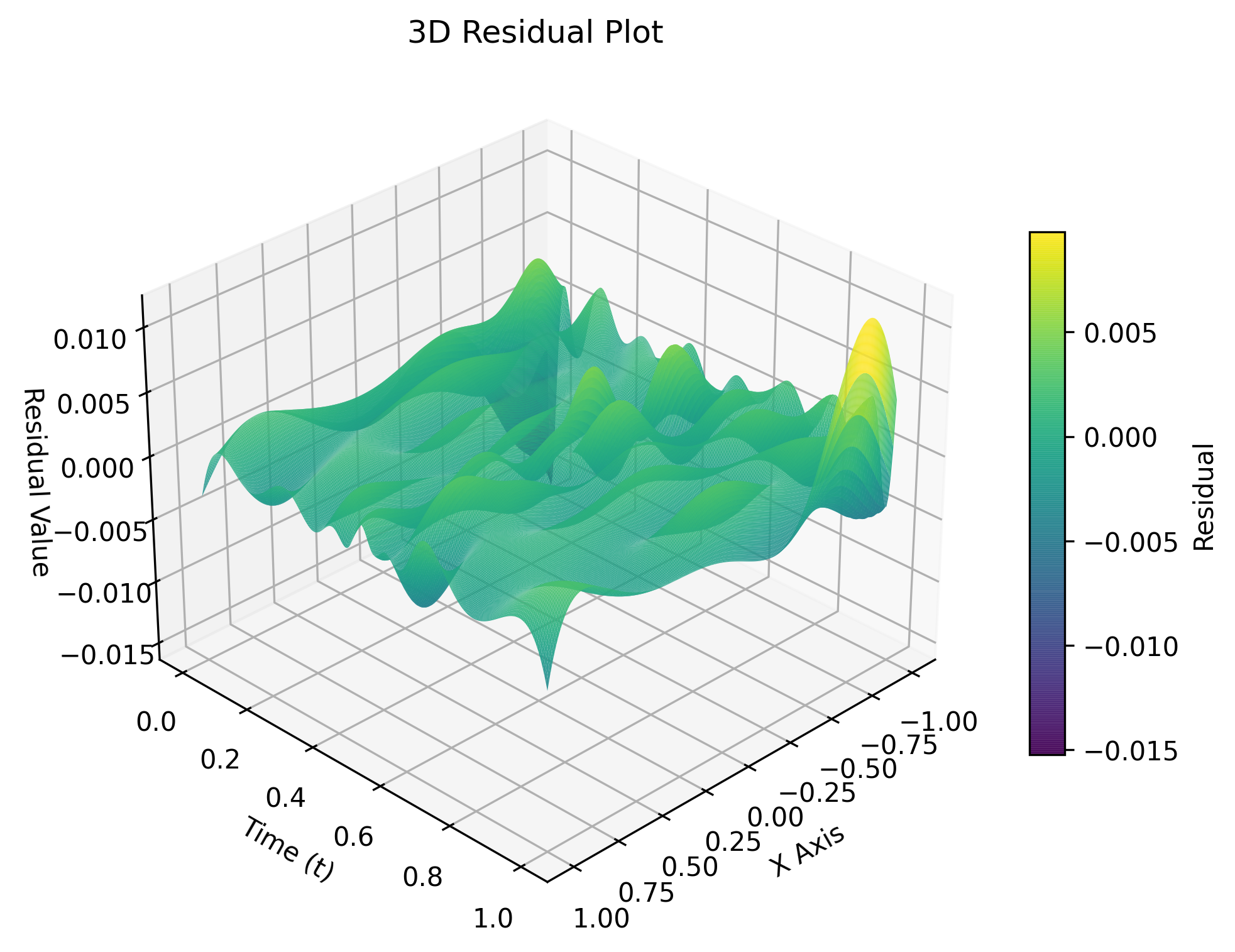}
\end{subfigure}
\hfill
\begin{subfigure}{0.32\textwidth}
    \centering
    \includegraphics[width=\textwidth]{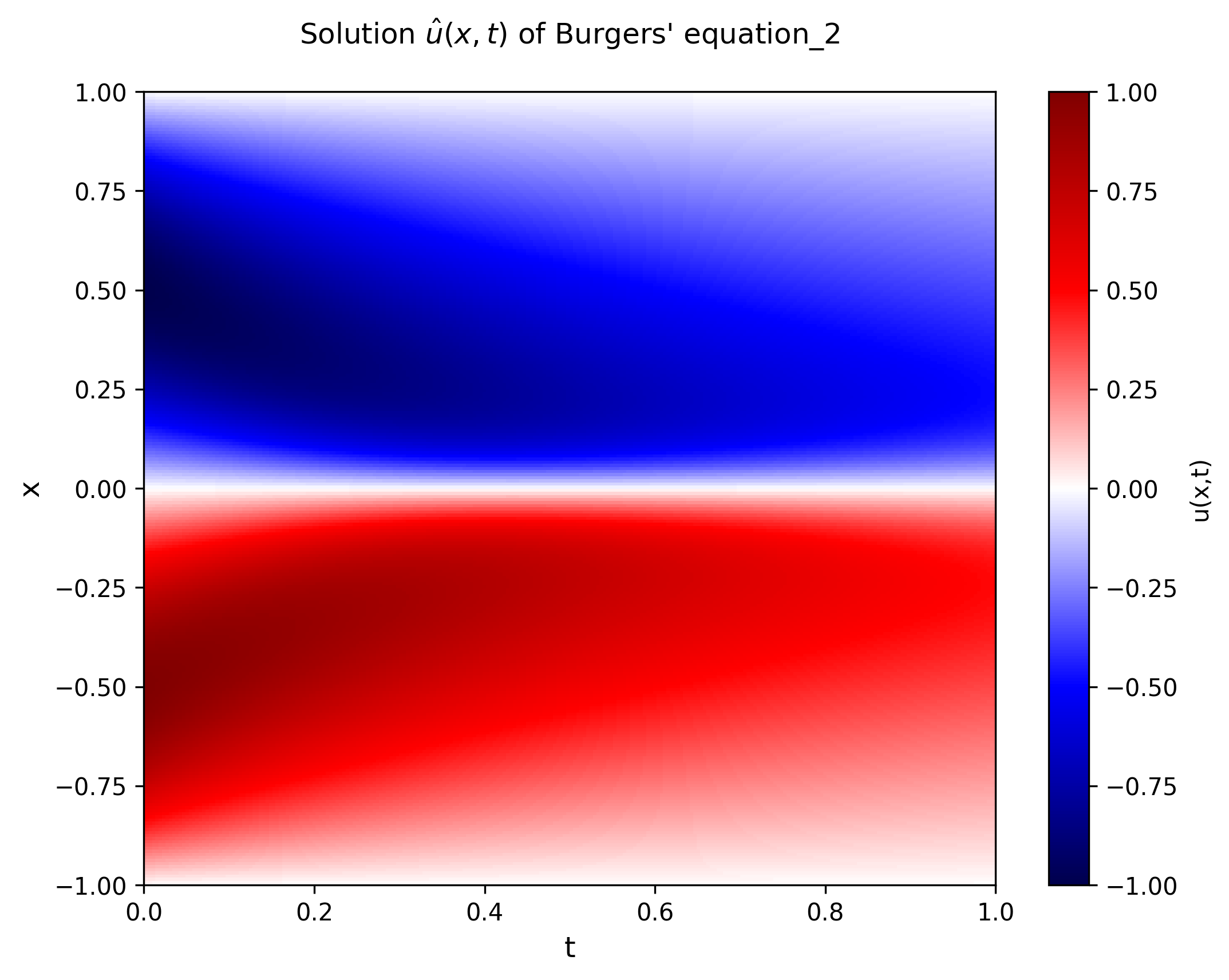}
\end{subfigure}
\end{minipage}

\vspace{1em} 

\begin{minipage}{\textwidth}
\centering
\footnotesize (c) Third-stage

\begin{subfigure}{0.32\textwidth}
    \centering
    \includegraphics[width=\textwidth]{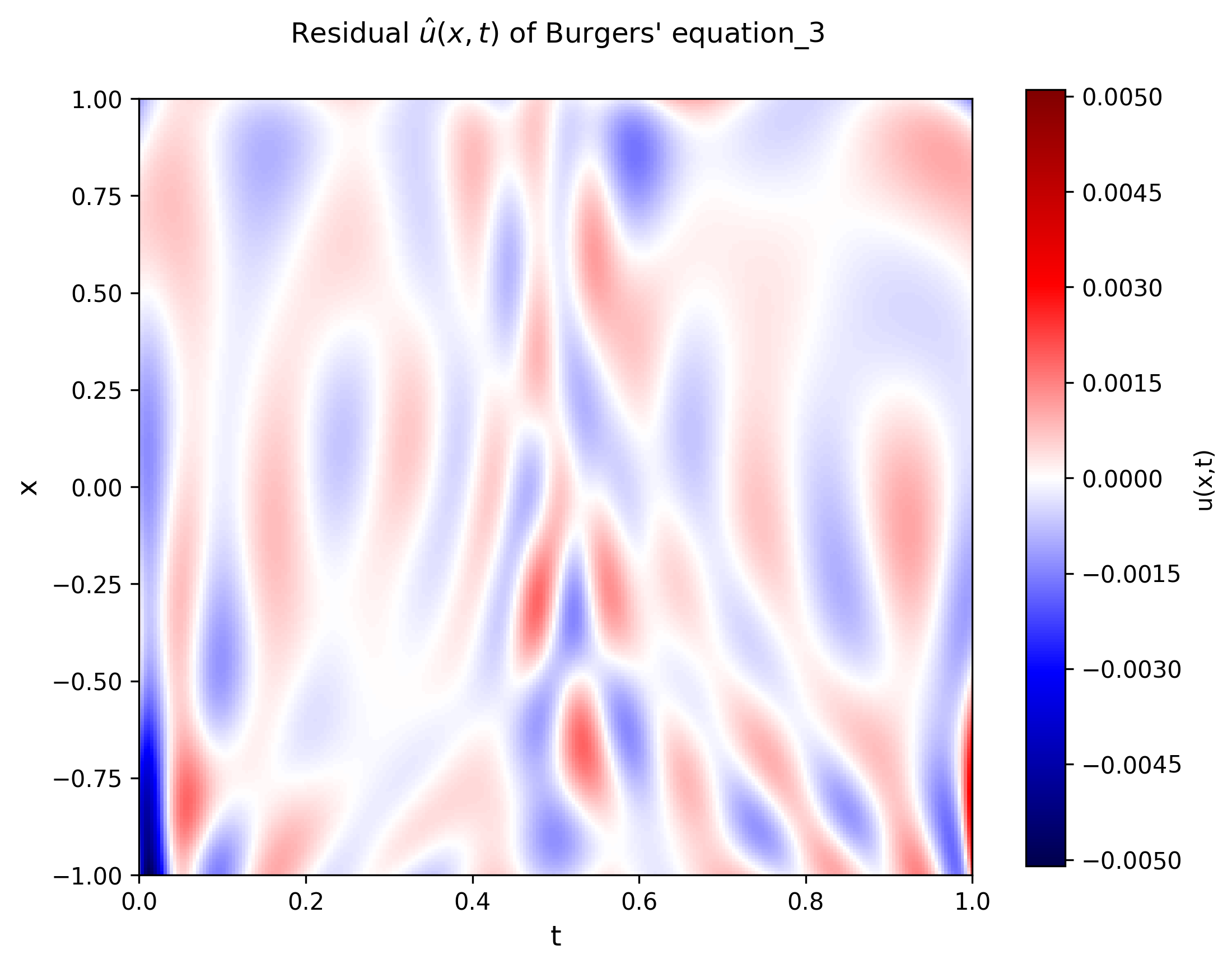}
\end{subfigure}
\hfill
\begin{subfigure}{0.32\textwidth}
    \centering
    \includegraphics[width=\textwidth]{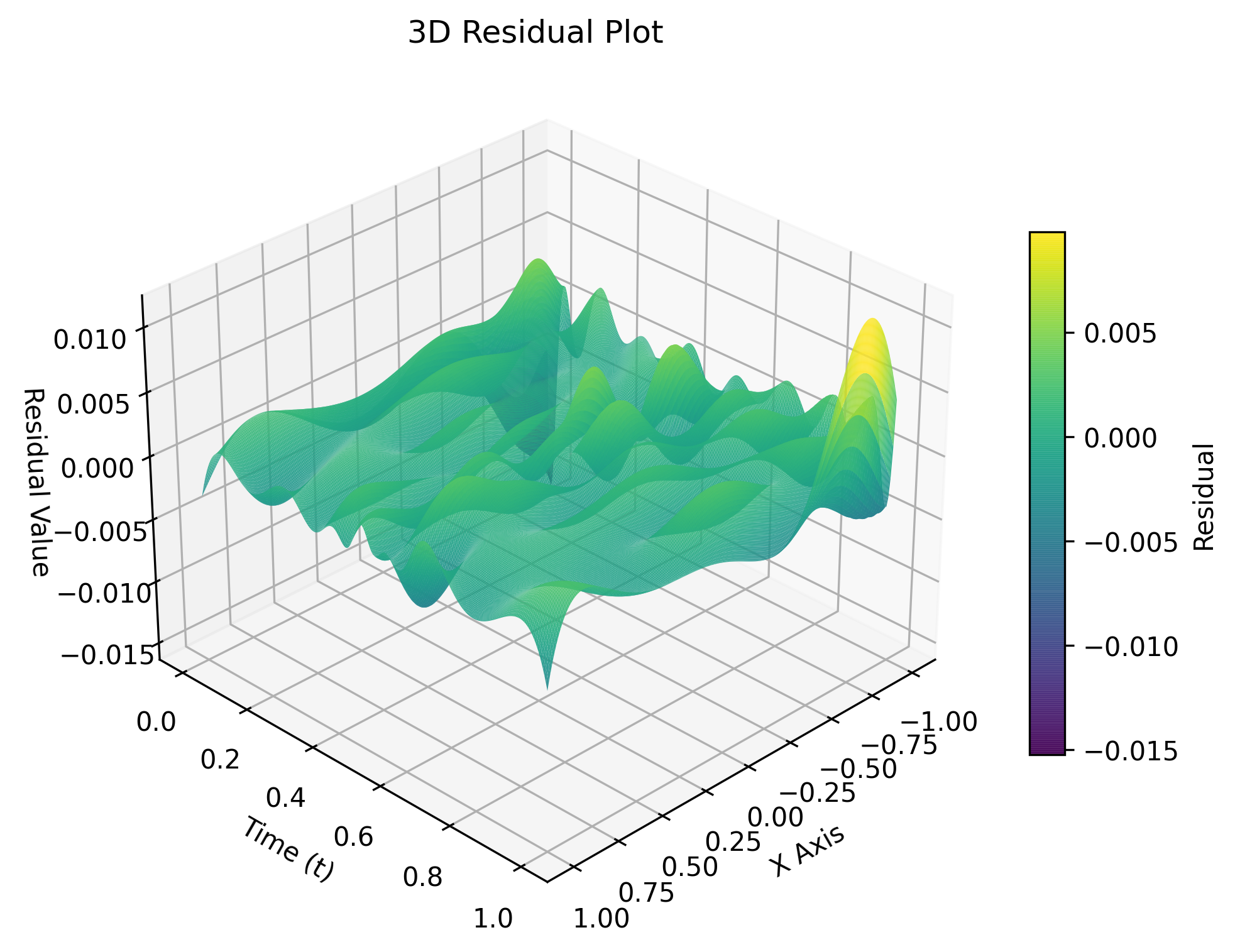}
\end{subfigure}
\hfill
\begin{subfigure}{0.32\textwidth}
    \centering
    \includegraphics[width=\textwidth]{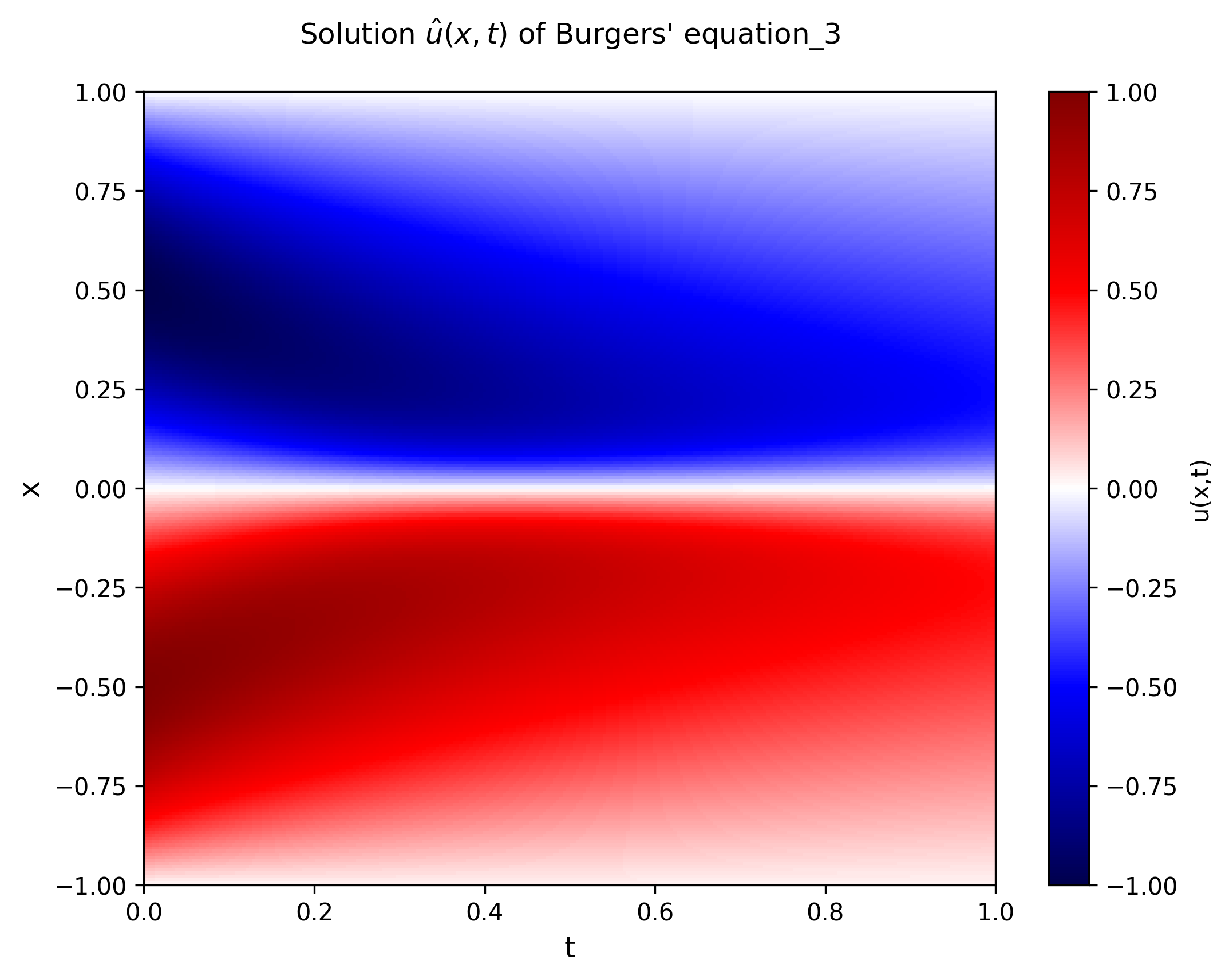}
\end{subfigure}
\end{minipage}
\caption{SI-MSPINNs for Burgers equation. (The left and middle figures show the two-dimensional and three-dimensional residual results of each stage, while the right figure shows the solution of the Burgers equation.)}
\label{fig:SI-MSPINNs-Burgers}
\end{figure}

The results of SI-MSPINNs solving the Burgers equation are shown in Figures \ref{fig:SI-MSPINNs-Burgers}.

\subsubsection{Results of the RFF-MSPINNs}
\begin{figure}[H]  
\centering
\begin{minipage}{\textwidth}
\centering
\footnotesize (a) First-stage

\begin{subfigure}{0.32\textwidth}
    \centering
    \includegraphics[width=\textwidth]{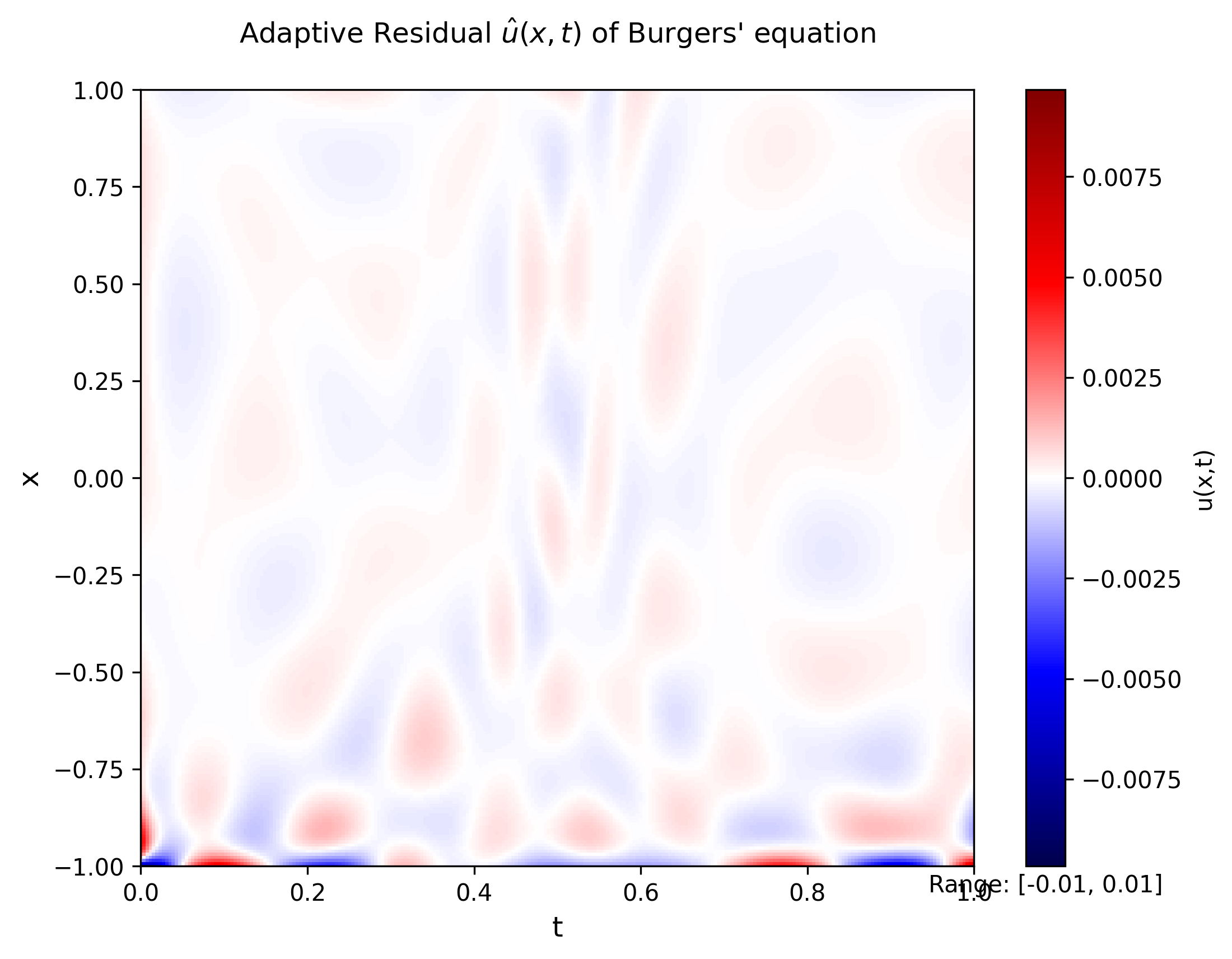}
\end{subfigure}
\hfill
\begin{subfigure}{0.32\textwidth}
    \centering
    \includegraphics[width=\textwidth]{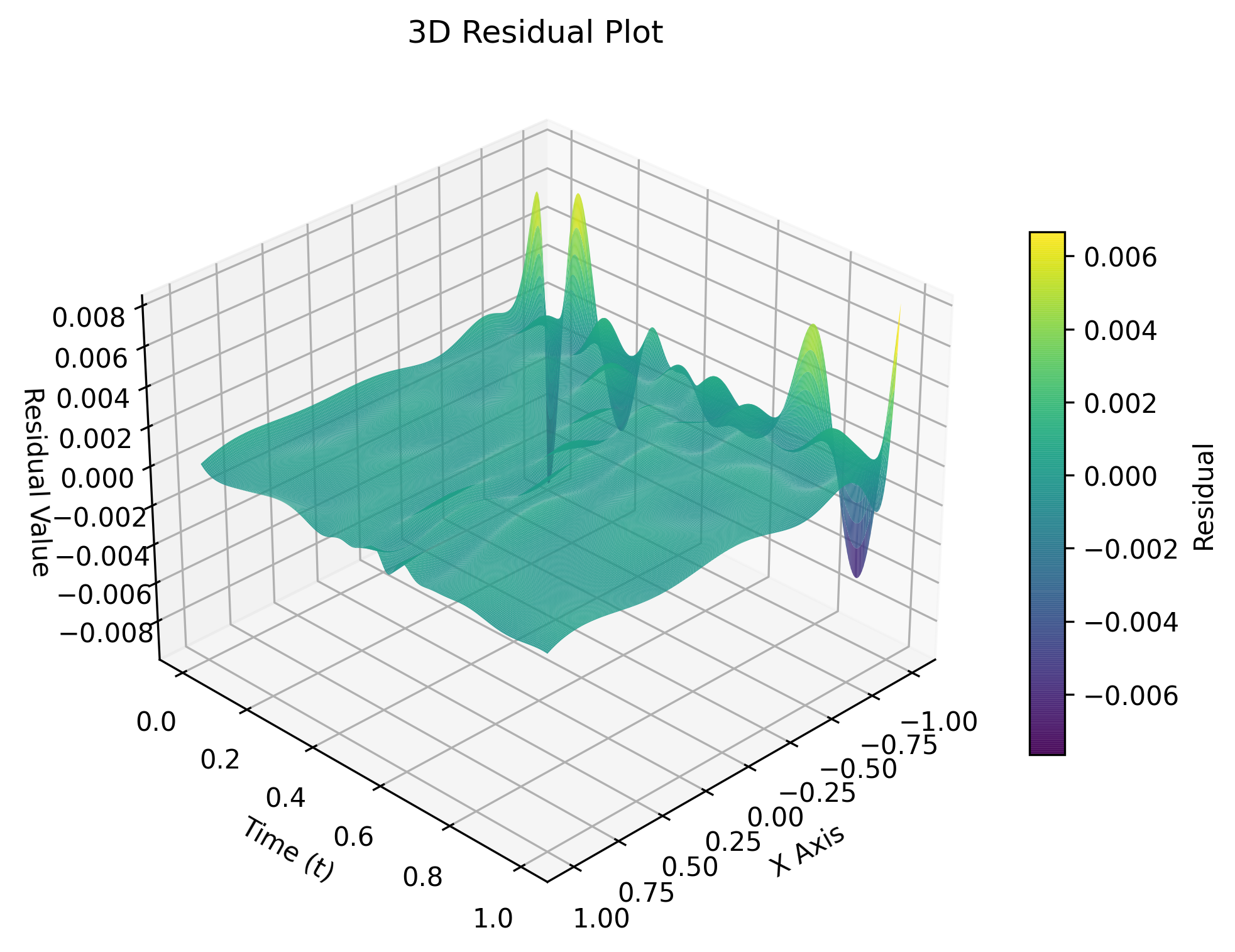}
\end{subfigure}
\hfill
\begin{subfigure}{0.32\textwidth}
    \centering
    \includegraphics[width=\textwidth]{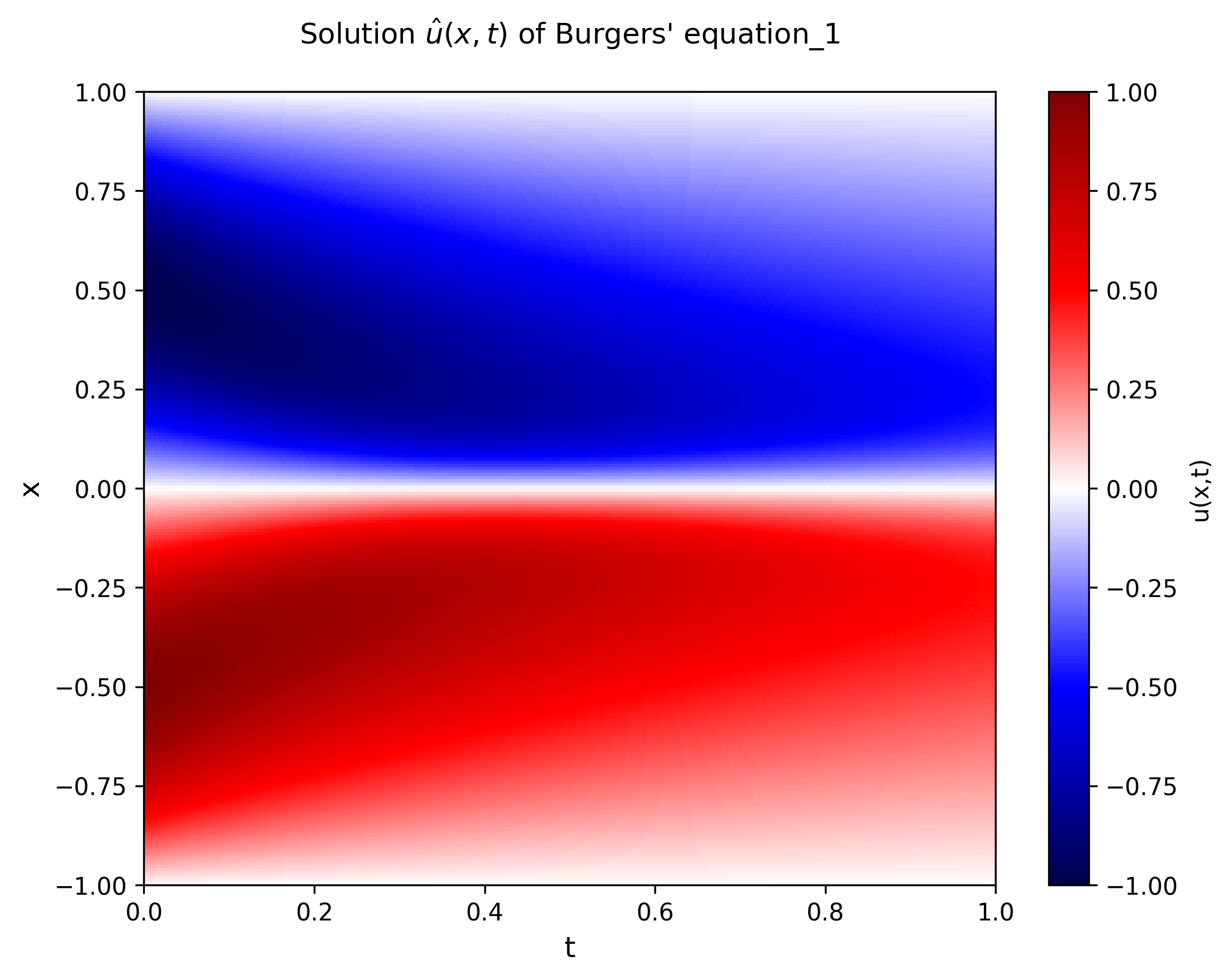}
\end{subfigure}
\end{minipage}

\vspace{1em} 

\begin{minipage}{\textwidth}
\centering
\footnotesize (b) Second-stage

\begin{subfigure}{0.32\textwidth}
    \centering
    \includegraphics[width=\textwidth]{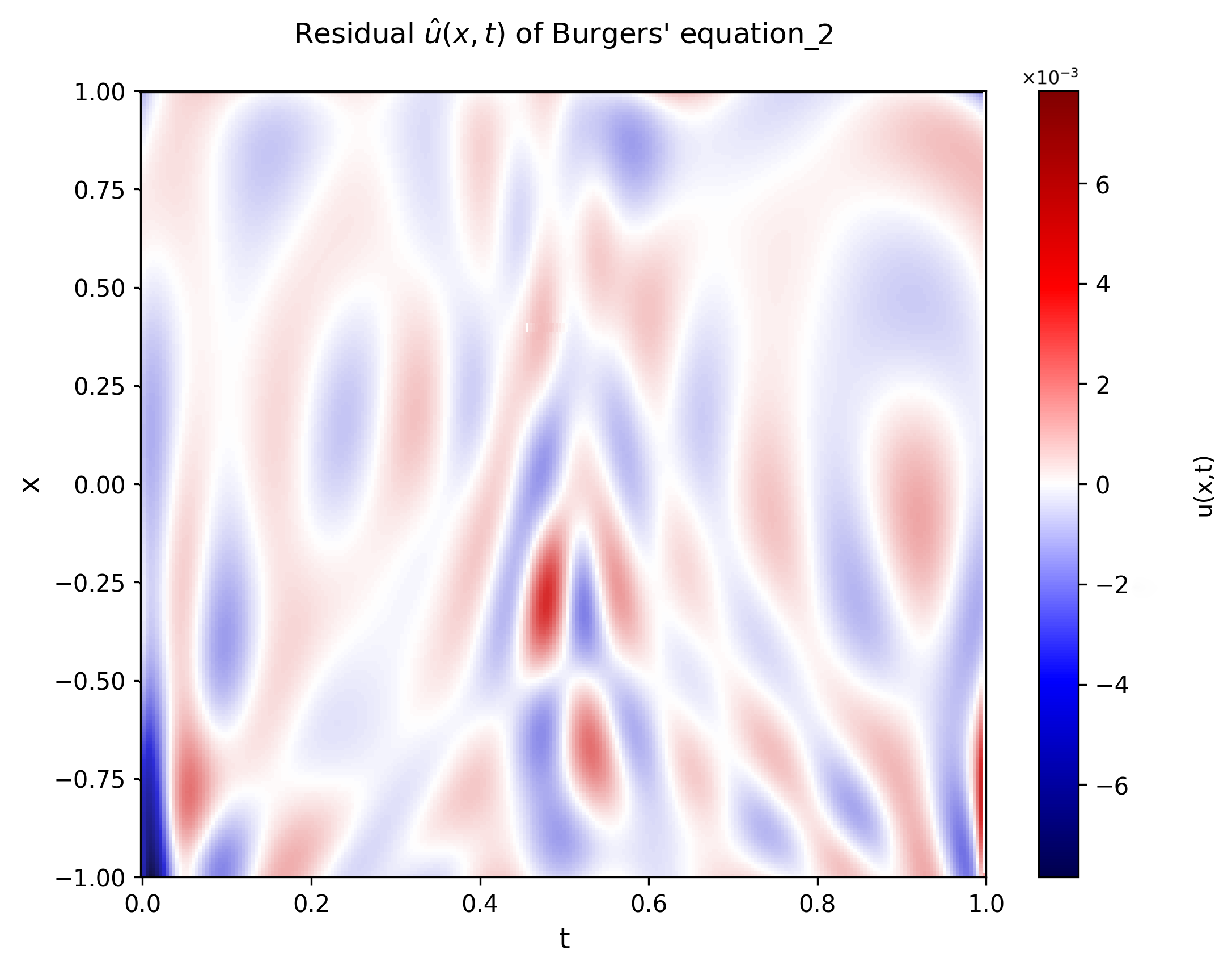}
\end{subfigure}
\hfill
\begin{subfigure}{0.32\textwidth}
    \centering
    \includegraphics[width=\textwidth]{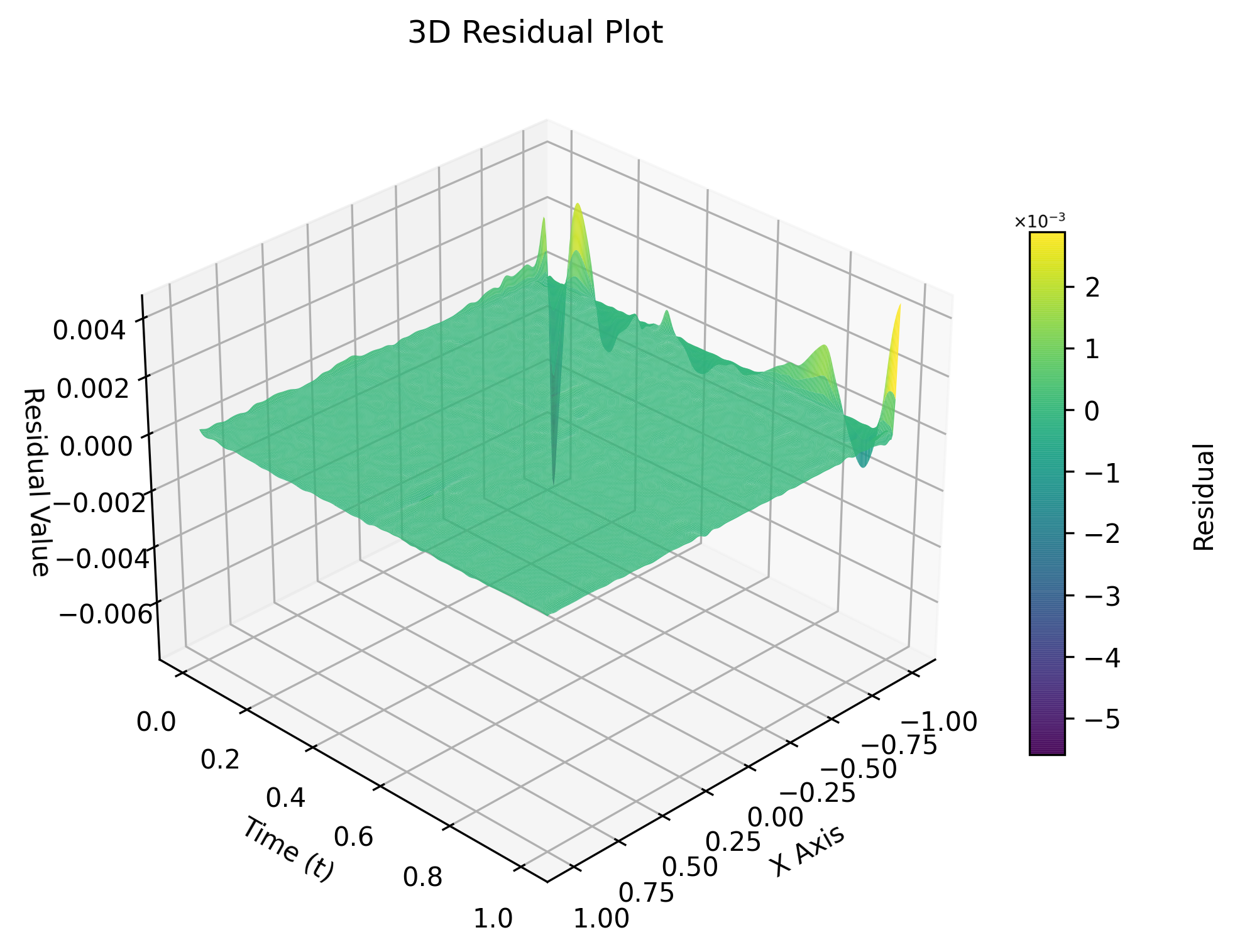}
\end{subfigure}
\hfill
\begin{subfigure}{0.32\textwidth}
    \centering
    \includegraphics[width=\textwidth]{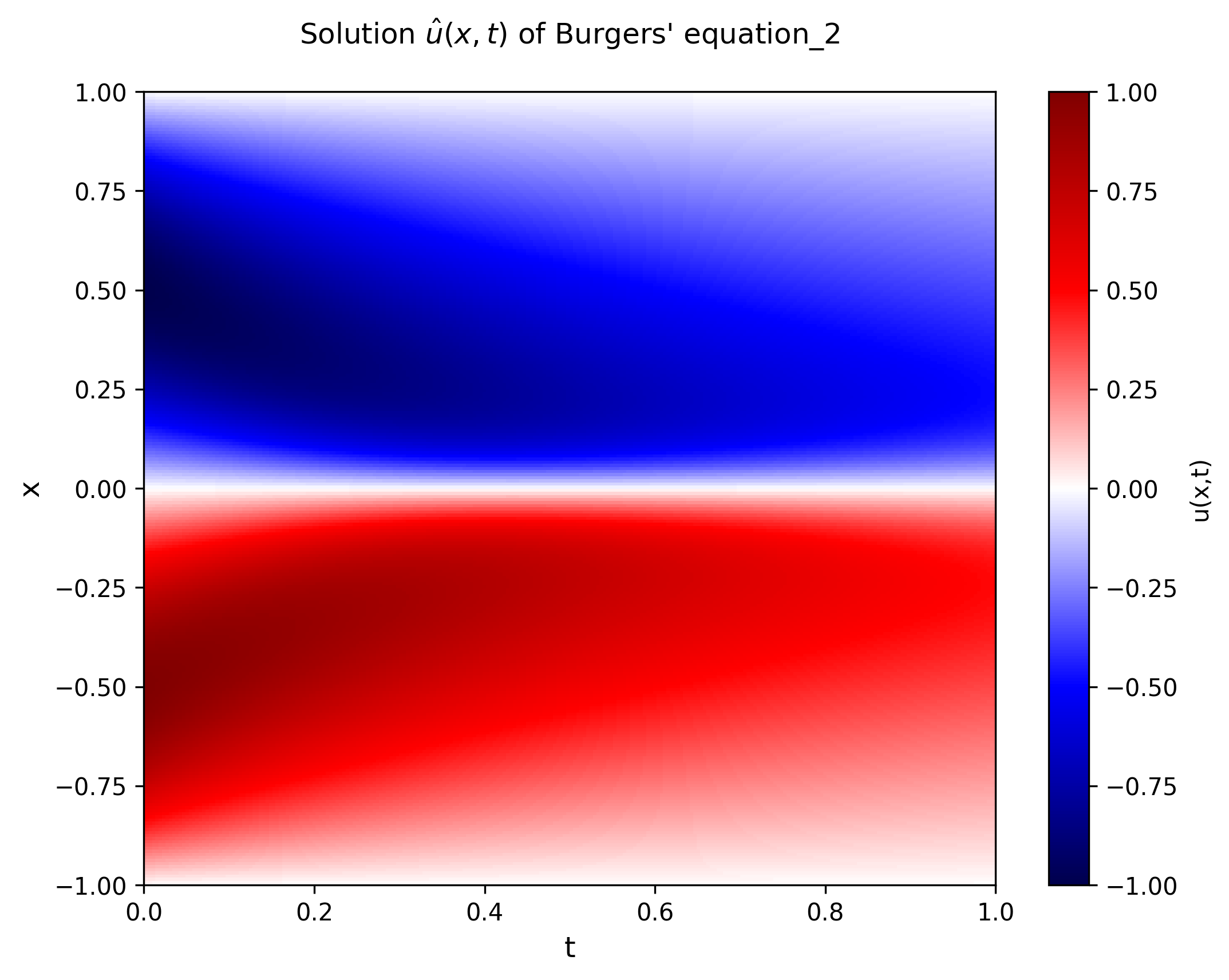}
\end{subfigure}
\end{minipage}

\vspace{1em} 

\begin{minipage}{\textwidth}
\centering
\footnotesize (c) Third-stage

\begin{subfigure}{0.32\textwidth}
    \centering
    \includegraphics[width=\textwidth]{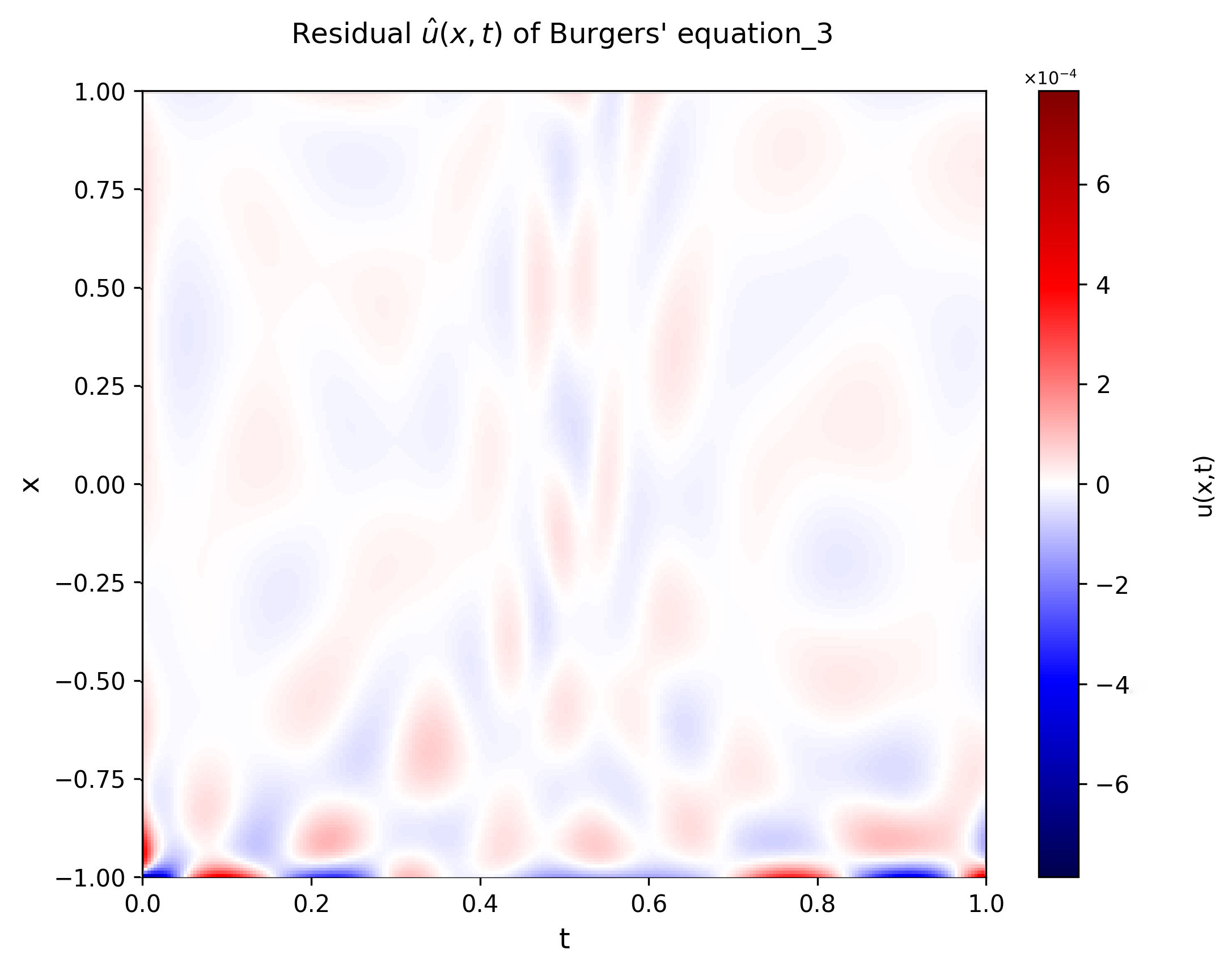}
\end{subfigure}
\hfill
\begin{subfigure}{0.32\textwidth}
    \centering
    \includegraphics[width=\textwidth]{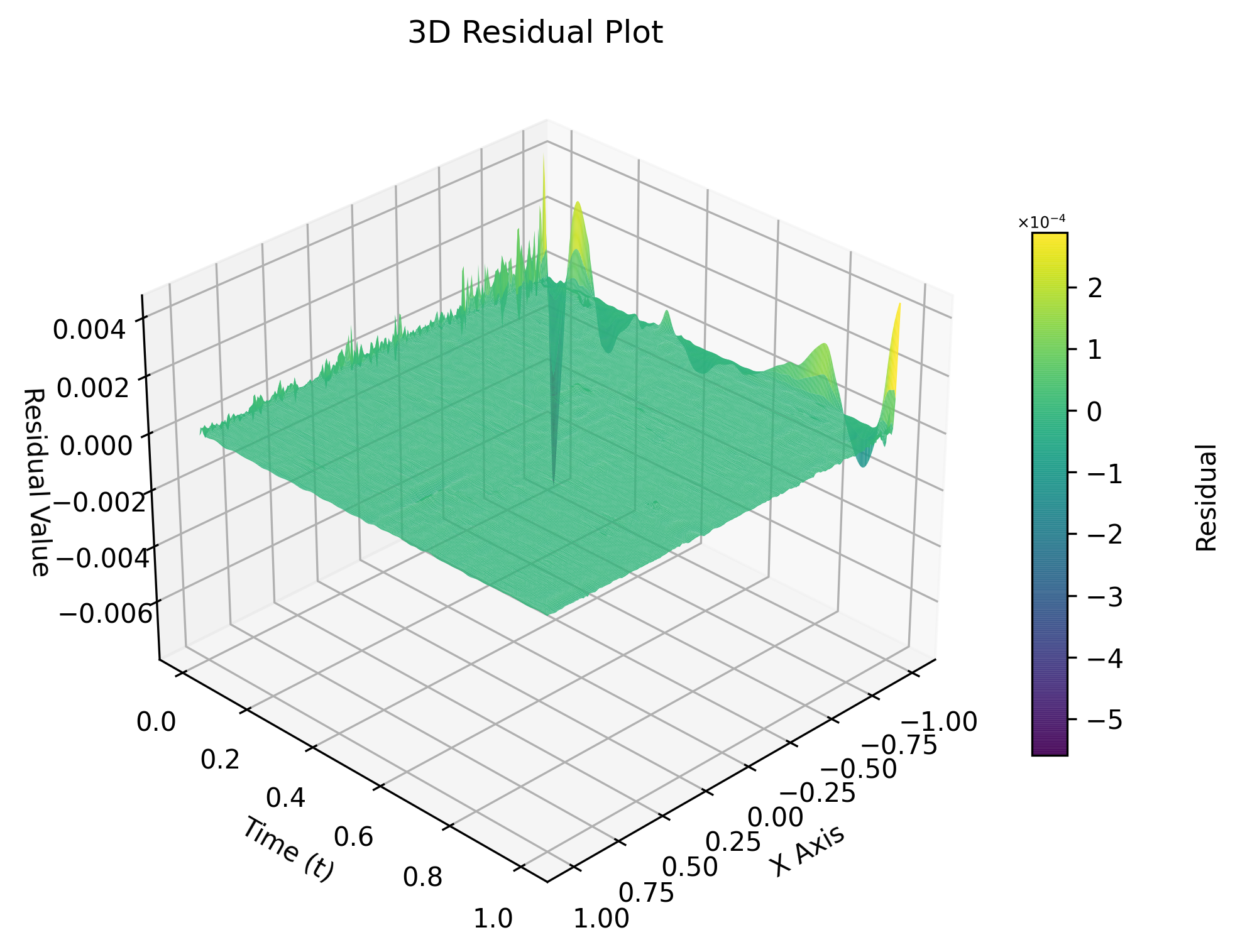}
\end{subfigure}
\hfill
\begin{subfigure}{0.32\textwidth}
    \centering
    \includegraphics[width=\textwidth]{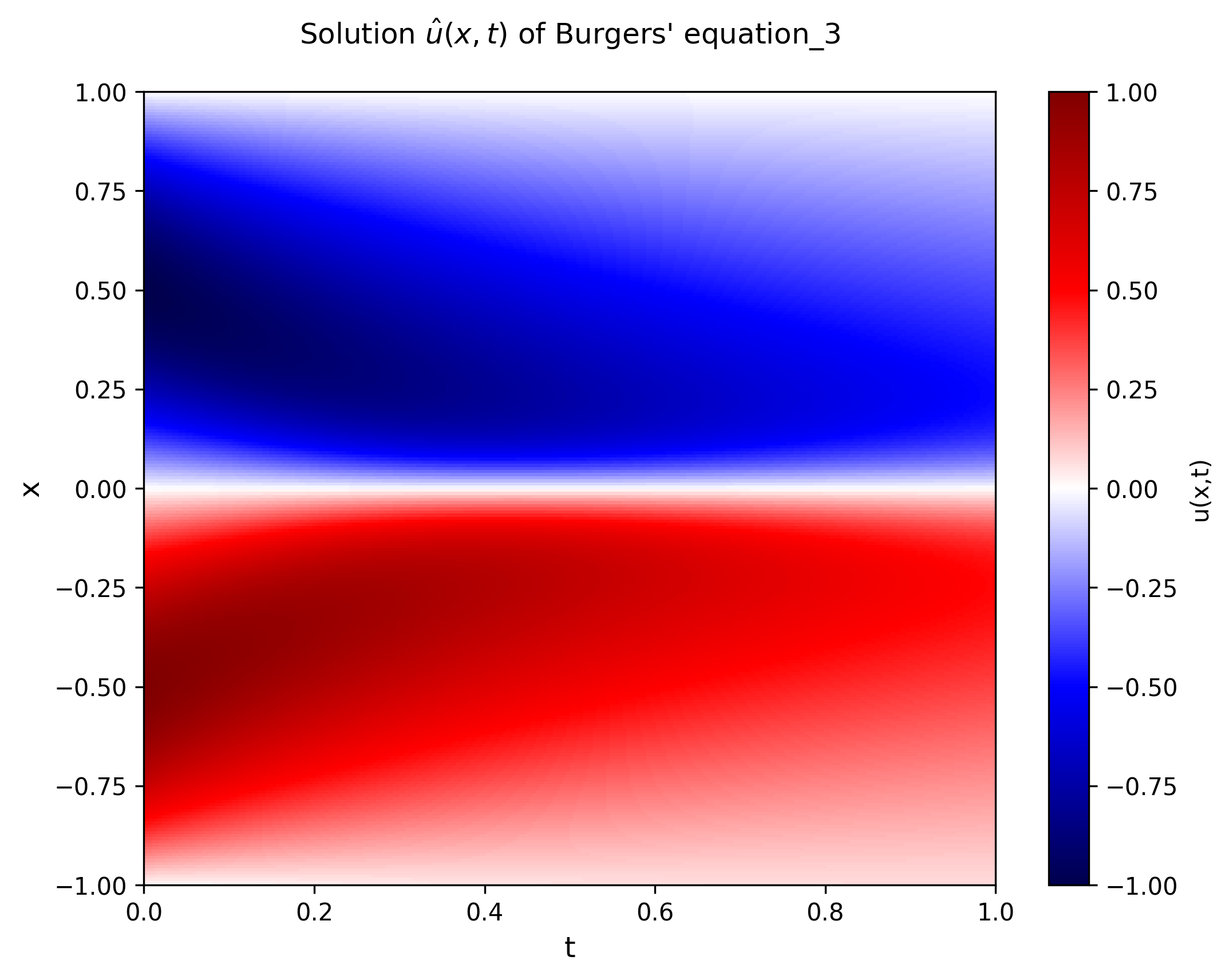}
\end{subfigure}
\end{minipage}
\caption{RFF-MSPINNs for Burgers equation. (The left and middle figures show the two-dimensional and three-dimensional residual results of each of the three stages, while the right figure shows the solution of the Burgers equation.)}
\label{fig:RFF-MSPINNs-Burgers}
\end{figure}

The results of RFF-MSPINNs solving the Burgers equation are shown in Figures \ref{fig:RFF-MSPINNs-Burgers}.


From Figures \ref{fig:SI-MSPINNs-Burgers} and \ref{fig:RFF-MSPINNs-Burgers}, we can observe that both the RFF-MSPINNs and SI-MSPINNs can improve accuracy through multistage residual learning, and enhance high-frequency capture capability by guiding multistage PINNs through spectral priors. In addition, RFF-MSPINNs does not directly use frequency components to construct a deterministic initialization in the first layer of the neural network, but uses frequency to guide the training of the neural network. Therefore, the final result of the third stage of RFF-MSPINNs has higher accuracy than SI-MSPINNs by an order of magnitude.

\begin{figure}[H]
\centering
\begin{subfigure}{0.49\textwidth}  
    \centering
    \includegraphics[scale=0.35]{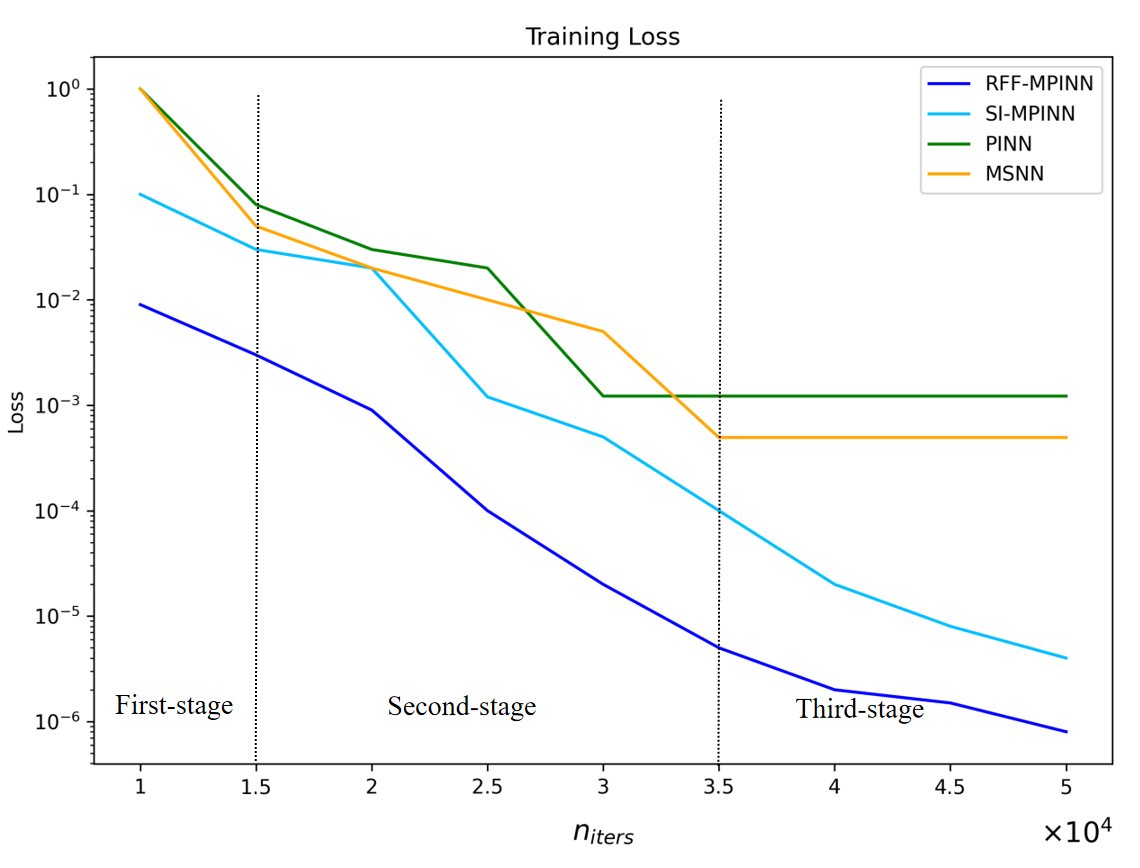}  
    \caption{Training loss figure for different methods.}
    \label{fig:burgers_training_loss}
\end{subfigure}
\hfill
\begin{subfigure}{0.49\textwidth}  
    \centering
    \includegraphics[scale=0.35]{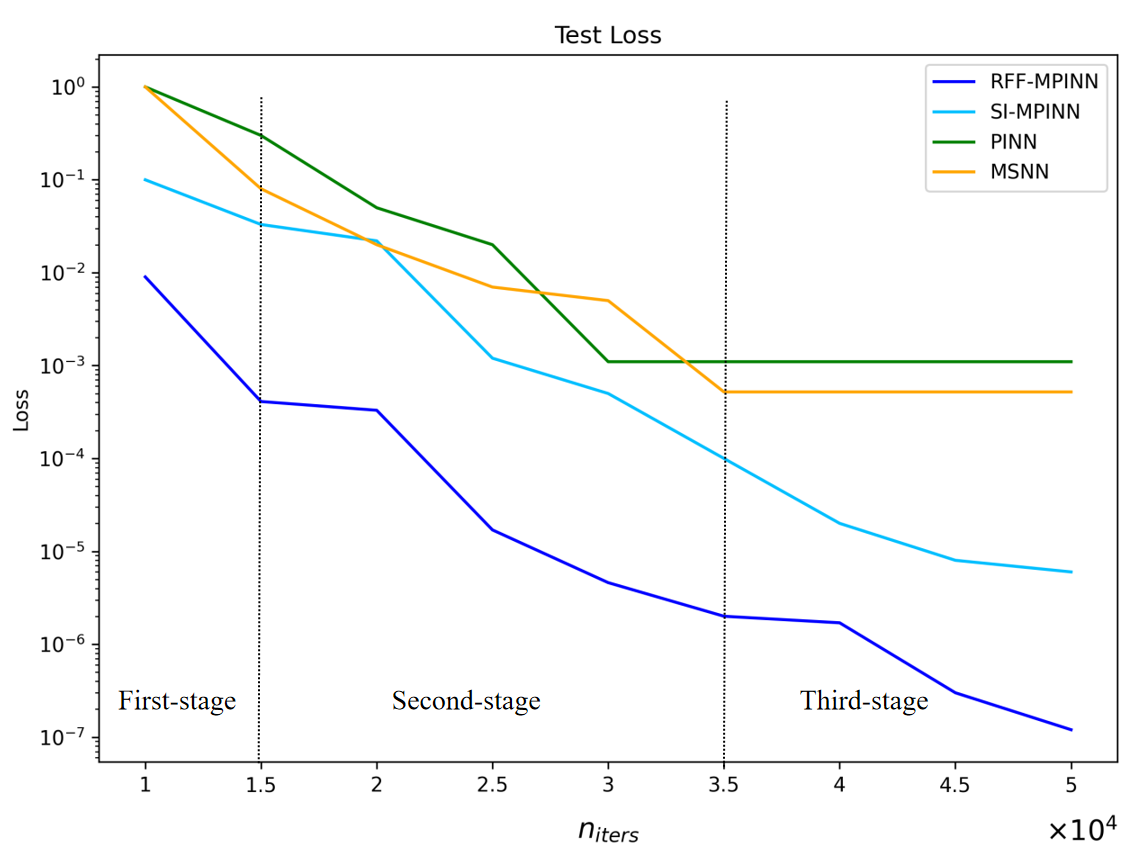}
    \caption{Test loss figure for different methods.}
    \label{fig:burgers_test_loss}
\end{subfigure}

\caption{Comparison figure of training loss and test loss for different methods.}  
\label{fig:burgers-losses}  
\end{figure}

As shown in Figures \ref{fig:burgers_training_loss} and \ref{fig:burgers_test_loss}, the loss values of traditional PINN and MSNN are $1.22 \times 10^{-3}$ and $4.94 \times 10^{-4}$, respectively. Meanwhile, the training loss values of SI-MSPINNs and RFF-MSPINNs in the third stage can reach $6.91 \times 10^{-6}$ and $7.34 \times 10^{-7}$, respectively. The three stages of training for SI-MSPINNs and RFF-MSPINNs took a total of 1326 seconds and 711 seconds, respectively.

The RFF-MSPINNs guides the network to focus on important frequencies while maintaining randomness, and is a method that falls between complete certainty and complete randomness. It does not require sorting all frequencies, but only sampling based on probability distribution. So the convergence speed of RFF-MSPINNs is faster.

\subsection{Helmholtz equation}
Consider the electromagnetic scattering problem of a dielectric disk described by the two-dimensional Helmholtz equation, defined as follows:
\begin{equation}
\left\{\begin{array}{l}
\frac{\partial^2 E_{r z}}{\partial x^2}+\frac{\partial^2 E_{r z}}{\partial y^2}+\omega^2 \mu \epsilon E_{r z}=0 \\
\frac{\partial^2 E_{i z}}{\partial x^2}+\frac{\partial^2 E_{i z}}{\partial y^2}+\omega^2 \mu \epsilon E_{i z}=0
\end{array}\right.,
\label{equ:helmholtz}
\end{equation}
where $E_{r z}$ and $E_{i z}$ represent the real and imaginary parts of the electric field intensity $E$ in the $z$-direction, respectively. $\epsilon$ and $\mu$ are the permittivity and permeability of the material, satisfying $\epsilon =\epsilon_0 \epsilon_r$ and $\mu =\mu_0 \mu_r$, where $\epsilon_0$ and $\mu_0$ are the permittivity and permeability of vacuum, $\epsilon_r$ and $\mu_r$ are the relative permittivity and relative permeability of the material. 

Since this scattering problem is defined in open domain, we introduce an artificial boundary to truncate it. On this boundary, a first-order absorbing boundary condition is imposed:
\begin{equation}
   \left\{\begin{array}{l}
\mathbf{n} \times \mathbf{E}_r-\sqrt{\frac{\mu}{\epsilon}} \mathbf{n} \times\left(\mathbf{H}_r \times \mathbf{n}\right)=\mathbf{n} \times \mathbf{E}_{r, i n c}-\sqrt{\frac{\mu}{\epsilon}} \mathbf{n} \times\left(\mathbf{H}_{r, i n c} \times \mathbf{n}\right) \\
\mathbf{n} \times \mathbf{E}_i-\sqrt{\frac{\mu}{\epsilon}} \mathbf{n} \times\left(\mathbf{H}_i \times \mathbf{n}\right)=\mathbf{n} \times \mathbf{E}_{i, \text { inc }}-\sqrt{\frac{\mu}{\epsilon}} \mathbf{n} \times\left(\mathbf{H}_{i, i n c} \times \mathbf{n}\right)
\end{array}\right. ,
\label{equ:abc}
\end{equation}
where $\mathbf{n}$ is the unit outward normal vector on the truncated boundary, $\mathbf{H}_r$ and $\mathbf{H}_i$ represent the real and imaginary parts of the magnetic field $\mathbf{H}$, $\mathbf{E^{inc}}$ and $\mathbf{H^{inc}}$ represent the incident electric and magnetic fields, respectively.

One can obtain the analytical solution of this problem \cite{chew1995waves}:
\begin{equation}
    \left\{  
\begin{aligned}  
u(r,\theta) = \sum_{n=-\infty}^{+\infty}  
i^n \bigl[J_n(Kr) + \alpha_n H_n^{(1)}(Kr)\bigr] \, e^{i n \theta}  && \quad \text{for } r \le R\\
u(r,\theta)  
= \sum_{n=-\infty}^{+\infty}  
   i^n \,\beta_n \, J_n\bigl(Kr\bigr)\, e^{i n \theta}  && \quad \text{for } r \ge R  
\end{aligned}  ,
\right.  
\label{equ18}
\end{equation}
where $\alpha_n$ and $\beta_n$ have the following forms:
\begin{equation}
   \alpha_n   
= \frac{  
   \mu\,K_i\,J_n'(K_i R)\,J_n(KR)  
   \;-\;  
   K\,J_n(K_i R)\,J_n'(KR)  
}{  
   K\,H_n^{(1)'}(KR)\,J_n(K_i R)  
   \;-\;  
   \mu\,K_i\,J_n'(K_i R)\,H_n^{(1)}(KR)  
} , 
\label{equ19}
\end{equation}
\begin{equation}
    \beta_n   
= \frac{  
   K\,H_n^{(1)'}(KR)\,J_n(KR)  
   \;-\;  
   K\,H_n^{(1)}(KR)\,J_n'(KR)  
}{  
   K\,H_n^{(1)'}(KR)\,J_n(K_i R)  
   \;-\;  
   \mu\,K_i\,J_n'(K_i R)\,H_n^{(1)}(KR)  
} , 
\label{equ20}
\end{equation}
where $R$ is the radius of the disk, $H_n^{(1)}$ and $H_n^{(1)'}$ are the first-kind Hankel function and its derivative, $J_n$ and $J_{n}'$ are the first-kind Bessel function and its derivative. This analytical solution is used as the ground truth for evaluating the accuracy of the methods.

\begin{figure}[htbp]
\centering
\includegraphics[scale=1]{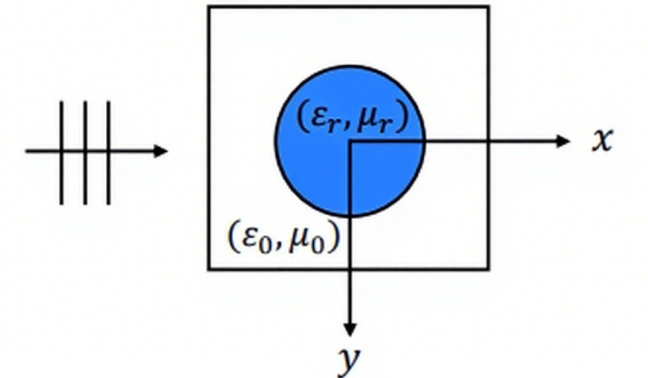}
\caption{Scattering of plane waves by a dielectric disk}
\label{disk}
\end{figure}

Using a plane wave as the excitation source $ H_{x}^{\mathrm{inc}}(x,y) = 0,H_{y}^{\mathrm{inc}}(x,y) =\cos\bigl(kx\bigr),E_{z}^{\mathrm{inc}}(x,y) = -\cos\bigl( kx\bigr)$ with a frequency $f=300\text{MHz}$, where $k=\frac{\omega}{c}$ is the wavenumber, $c$ is the speed of light in vacuum, $\omega=2\pi f$ is the angular frequency, and $f$ is the wave frequency. The computational domain is $\Omega=[-1,1]^2$, and the disk radius is 0.25, as shown in Figure \ref{disk}. This equation is difficult to solve for PINNs method, because we could have discontinuity of the boundary of the disk.

Set the number of training residual point samples within the domain to 2540, and the number of boundary sampling points to 80. We still use a three-stage architecture. All the three networks use fully connected neural networks with a depth of 4 (equivalent to 3 hidden layers) and a width of 20, and the optimizer is Adam+L-BFGS for all the three stages, the same as those in the previous example. 

\subsubsection{Results of the SI-MSPINNs}
\begin{figure}[htbp]
\centering
\includegraphics[scale=0.25]{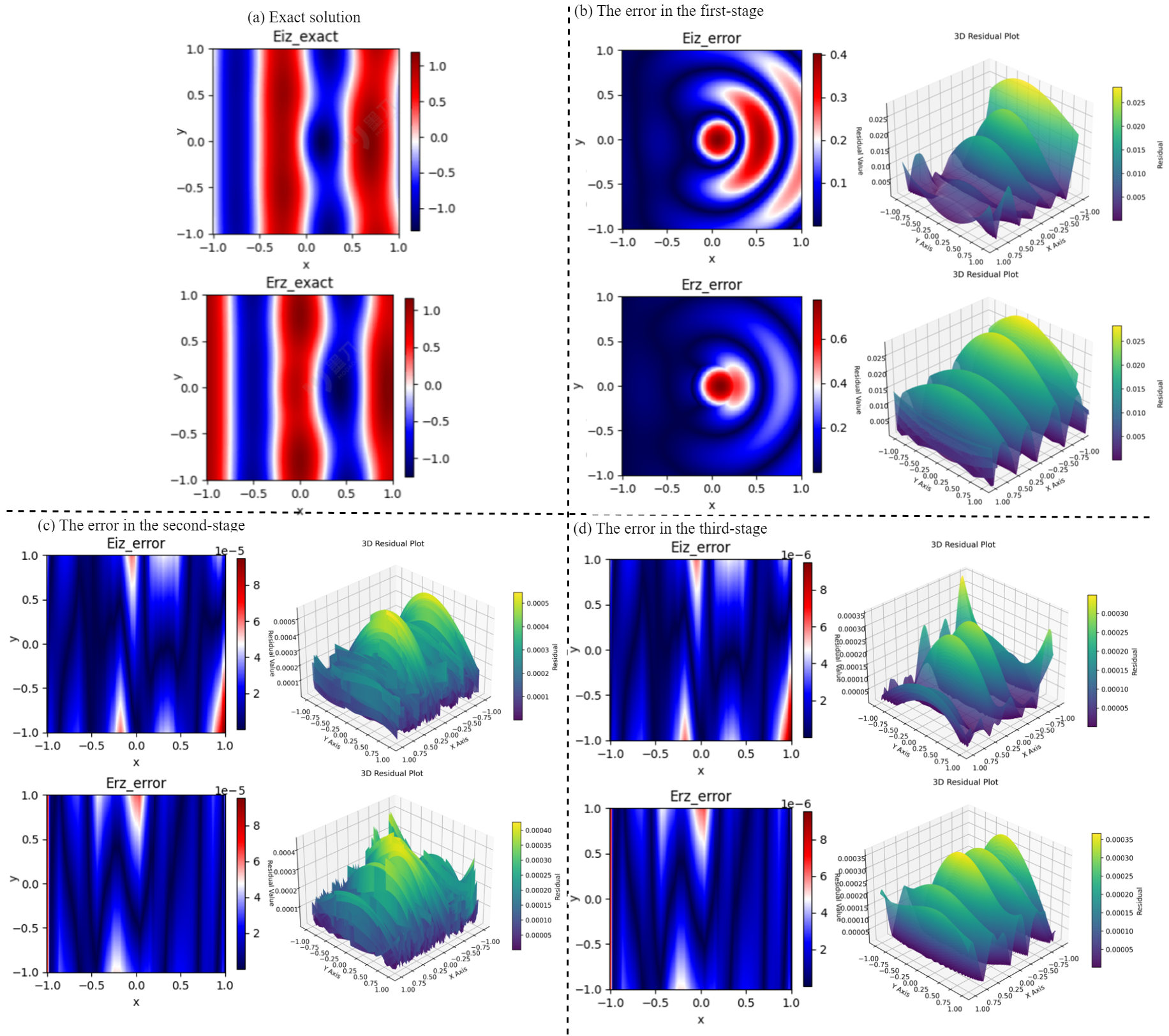}
\caption{SI-MSPINNs for Helmholtz equation ($\epsilon=1$).}
\label{h1}
\end{figure}

When $\epsilon=1$, the results of SI-MSPINNs are shown in Figure \ref{h1}.

\subsubsection{Results of the RFF-MSPINNs}
\begin{figure}[htbp]
\centering
\includegraphics[scale=0.25]{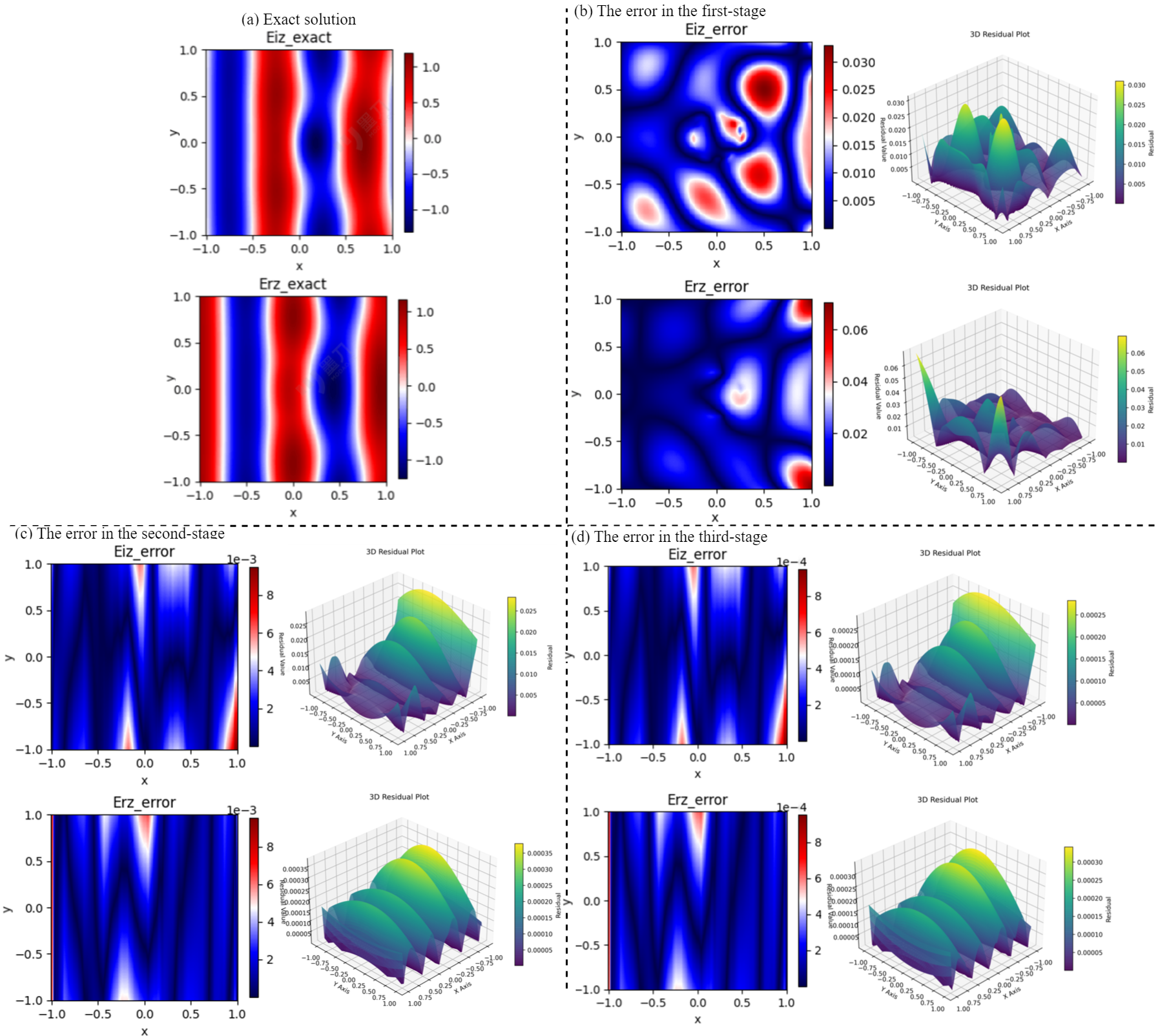}
\caption{SI-MSPINNs for Helmholtz equation($\epsilon=1.5$).}
\label{h2}
\end{figure}

In Figure \ref{h2} we show the results of RFF-MSPINNs when $\epsilon=1.5$.

In Table \ref{T2}, we compare the $L_2$ error of PINN, MSNN, SI-MSPINNs, and RFF-MSPINNs. We can see that both SI-MSPINNs and RFF-MSPINNs effectively improve the accuracy of the original PINNs as well as MSNN. On can also observe that RFF-MSPINNs has higher accuracy, because the frequency sampling of RFF-MSPINNs is proportional to the PSD of the objective function, prioritizing the coverage of critical frequencies. This ensures that more features are allocated to important frequencies and fewer to low-energy frequencies. In this way, the neural network is initialized with features that match the spectrum of the objective function, enabling it to learn the main frequency components of the objective function more quickly, especially the high-frequency components with higher energy. Furthermore, through multistage training, residuals are gradually refined, with each stage initialized based on the spectral information of residuals from the previous stage. This collaborative design allows RFF-MSPINNs to more effectively capture both low-frequency and high-frequency features, address interface discontinuity problems, and thus achieve high accuracy.

\begin{table}[htbp]  
    \centering  
    \caption{Comparison of results of different methods.}  
    \small  
    \begin{tabular}{lcccccc}  
        \toprule  
        method & \multicolumn{2}{c}{$\epsilon=1$} & \multicolumn{2}{c}{$\epsilon=1.5$} & \multicolumn{2}{c}{$\epsilon=2$} \\
        \cmidrule(lr){2-3} \cmidrule(lr){4-5}\cmidrule(lr){6-7}  
        & $L_2(E_{rz})$ & $L_2(E_{iz})$ & $L_2(E_{rz})$ & $L_2(E_{iz})$& $L_2(E_{rz})$ & $L_2(E_{iz})$ \\   
        \midrule  
        PINN       & $0.0092$ & $0.0084$ & $0.0206$ & $0.035$ & $0.207$ & $0.158$\\
        MSNN       & $0.0007$ & $0.0005$ & $0.0024$ & $0.0098$ & $0.034$ & $0.061$\\
        SI-MSPINNs   & $3.2 \times 10^{-6}$ & $2.1 \times 10^{-6}$ & $0.0019$ & $0.002$& $0.0075$ & $0.031$ \\
        RFF-MSPINNs  &$2.9 \times 10^{-6}$ & $4.1 \times 10^{-5}$&$1.7 \times 10^{-4}$ & $2.6 \times 10^{-4}$&$0.0012$ & $0.0003$ \\
        \bottomrule  
    \end{tabular}
    \label{T2}  
\end{table}

\section{Conclusion \& Future Work}\label{sec:Conclusion}
We propose the SI-MSPINNs and RFF-MSPINNs for solving PDEs with high accuracy. 
\begin{itemize}
    \item Compared with PINNs, both models effectively overcome their spectral bias issues. Specifically, SI-MSPINNs extracts dominant spectral patterns of residuals through discrete Fourier transform to guide network initialization, enabling accurate capture of high-frequency physical features. This allows the network to more efficiently adapt to high-frequency components when learning PDE residuals. RFF-MSPINNs, on the other hand, dynamically adjusts the frequency sampling distribution of random Fourier features based on residual power spectral density, prioritizing focus on high-energy physical modes. It breaks through the limitations of traditional PINNs, such as inefficient derivative calculation caused by automatic differentiation and insufficient ability to learn high-frequency features.
    \item Compared with MSNNs, SI-MSPINNs and RFF-MSPINNs further address their core defects in 2D problems. Although MSNNs alleviate spectral bias through multistage residual learning, their convergence efficiency decreases significantly in high-dimensional scenarios, and spectral initialization is susceptible to spectral aliasing, leading the network to fall into local optima. In contrast, SI-MSPINNs adopts a spectrum-informed hierarchical initialization strategy, enabling each stage of the network to specifically correct the residuals of the preceding cumulative solution. Combined with the gradual enhancement of the physical constraint weights, it achieves a progressive error reduction while ensuring physical consistency. RFF-MSPINNs replaces frequency sorting with probability sampling, avoiding the computational burden of sorting spectral features in high dimensions. It enables the network focus on key frequencies while maintaining randomness, resulting in faster convergence and stronger high-dimensional scalability. Results on solving Burgers equation and parameterized Helmholtz equation show the effectiveness of the proposed methods. 
\end{itemize}

The SI-MSPINNs and RFF-MSPINNs have shown higher accuracy in 2D problems. Our future work will be focused on breaking through the hierarchical spectrum initialization strategy for higher-dimenstion scenes, and improving the ability to capture local features in shock wave regions through the fusion of adaptive resampling \cite{RAR} and gradient enhancement constraints \cite{gPINN}. At the same time, exploring the use of wavelet basis functions instead of Fourier basis to optimize the spectral representation of non-stationary signals such as turbulence, and combining quantum Fourier transform to accelerate high-dimensional spectral processing.
\section{Acknowledgments}
The work is partially supported by Natural Science Foundation of Sichuan Province (Grant No. 24NSFSC1366) and NSFC (Grant No. 12101511).
\bibliographystyle{elsarticle-num}
\bibliography{reference}

\end{document}